\documentclass[12pt,a4paper]{amsart}
\usepackage{graphicx}
\usepackage{epstopdf}
\usepackage{amssymb}
\usepackage{amsthm}
\usepackage{amsmath}
\usepackage{pstcol, pst-node}
\usepackage{mathrsfs}
\usepackage{amscd}
\usepackage[all]{xy}
\newtheorem{thm}{Theorem}

\newtheorem{cor}[thm]{Corollary}
\newtheorem{prop}[thm]{Proposition}
\newtheorem{lem}[thm]{Lemma}
\newtheorem{rem}[thm]{Remark}

\theoremstyle{definition}
\newtheorem{defn}[thm]{Definition}

\newtheorem{prop-def}[thm]{Proposition-Definition}

\theoremstyle{remark}

\begin{document}
\title{Counting curves via degeneration}

\author{Takeo Nishinou}
\date{}
\thanks{email : nishinou@rikkyo.ac.jp}
\address{Department of Mathematics, Rikkyo University,
  Nishi-Ikebukuro, Toshima, Tokyo, Japan} 
\subjclass[2000]{}
\keywords{}
\maketitle
\begin{abstract}
We develop a technique to study curves in a variety which has a degeneration 
 into some union of varieties.
The class of such varieties is quite broad, but the theory 
 becomes particularly useful when the variety has a degeneration into 
 a union of toric varieties.
Hypersurfaces are typical examples, and we study lines 
 on quartic K3 surfaces and quintic Calabi-Yau
 hypersurfaces in detail.
In particular, we combinatorially prove the existence of 
 2875 lines in a generic quintic Calabi-Yau 3-fold.
Also, we give a geometric construction of walls in the Gross-Siebert construction 
 of Calabi-Yau varieties \cite{GS}.
\end{abstract}
\section{Introduction}
This paper is, in a sense, a sequel to the paper \cite{N}, where we studied 
 higher genus curves in a toric variety.
In this paper, we develop a technique which enable us to use ideas from 
 tropical geometry to study curves in a variety which is not necessarily toric.
Thus, combining with \cite{N}, we can theoretically study higher genus
 curves in such a variety by combinatorial way.
 
In \cite{N2}, we already used the theory of tropical curves to study curves in 
 varieties which degenerate to irreducible toric varieties.
Examples included flag varieties and some moduli space of bundles on 
 a Riemann surface.
In this paper, we consider the more general case
 that a variety has a degeneration to a \emph{union} of 
 toric varieties.
The class of such varieties is broad, for example, in mirror symmetry it is conjectured
 that any Calabi-Yau variety has such a degeneration. 
In fact, the method of this paper works even 
 when a variety has a degeneration to a union of 
 varieties which are not necessarily toric.
But we mainly concentrate
 on the case when the variety has a degeneration to a union of toric varieties
 ('toric degeneration' in the language of \cite{NS}), where the theory becomes
 particularly efficient.

Our principal examples are K3 and Calabi-Yau hypersurfaces.
In the next section, we study the K3 case, where the essential calculation
 of this paper is done.
The core of our calculation consists of 2 points:
\begin{itemize}
\item Calculate the log normal sheaf.
\item Geometrically understand the tangent and obstruction classes. 
\end{itemize}
The first point is rather simple calculation, and the method is easily  
 applicable to broad situations.
The second point is more complicated but interesting calculation.
Namely, it is the calculation of the deformation when the obstruction
 does not vanish (the Kuranishi map).
There is no general method to calculate the Kuranishi map, and 
 to perform the calculation is usually quite difficult in explicit examples.
We see that degeneration makes the calculation possible in some situations. 

As applications, we calculate the Kuranishi map for lines in 
 K3 surfaces, and 
 we also study lines in a quintic Calabi-Yau hypersurface.
There we borrow an idea of Sheldon Katz \cite{Katz}, which we
 learned from Mark Gross, giving 
 a way to combinatorially count
 2875 lines in it.
In Section \ref{sec:disks}, we give another example of our calculation, applied to 
 holomorphic disks.
Namely, we realize the walls in Gross
 and Siebert's construction as a family of holomorphic disks, 
 as conjectured by them \cite{GS}.
Other applications of the techniques in this paper can be found in 
 \cite{N, NY}.\\

\noindent
{\bf Acknowledgment.} 
The starting point of this paper was the combinatorial idea of counting lines 
 in a quintic Calabi-Yau 3-fold.
Bernd Siebert kindly showed me an email from Mark Gross, in which Gross 
 explained this idea.
On the other hand, Gross told me that it was an idea going back to Sheldon Katz
 \cite{Katz}.
We would like to thank all of them, without any of these coincidences this paper 
 might not have been written.
We would also like to thank Ionut Ciocan-Fontanine and
 for useful conversation.
The author was supported by JSPS KAKENHI Grant Number
 26400061.

\section{Lines on K3 surfaces via deformation}\label{sec:K3}
Consider a degeneration of a quartic K3 surface given by the equation
\[
\bar x\bar y\bar z\bar w+ tf = 0,
\]
where $\bar x, \bar y, \bar z, \bar w$ 
 are homogeneous coordinates of $\Bbb P^3$,
 $t\in\Bbb C$ is the parameter of the degeneration
 and $f$ is a generic homogeneous quartic polynomial
 of $\bar x, \bar y, \bar z, \bar w$.
Let 
\[
\mathfrak X\subset \Bbb P^3\times \Bbb C
\]
 be the variety defined by the above
 equation (the total space of the degeneration), and $X_0$ be the central fiber. 
Let 
\[
i_0\colon X_0\to\mathfrak X
\] 
 be the inclusion.
The space $X_0$ is the union of 4 projective planes glued along projective lines:
\[
X_0 = \cup_{i =1}^4{\Bbb P^2_i}
\]
Each $\Bbb P^2_i$ has a natural structure of a toric variety, in which the lines
 mentioned above are the toric divisors. 
Let 
\[
\ell_1, \dots, \ell_6
\]
 be these projective lines.
Let $L = \cup_{i=1}^6\ell_i$ be the union of them.
Each $\ell_i$ has 2 distinguished points, which are the triple-intersections
 of the projective planes.
We write by $\ell_i^{\circ}$ the complement of these 2 points.

Since $f$ is generic, the total space $\mathfrak X$ has 24 singular points, 
 and for each $i$, 4 of them lie on $\ell_i^{\circ}$.
Let 
\[
\mathcal S\subset X_0
\]
 be the set of these singular points.
These singular points have the same local structure:
\begin{lem}\label{lem:coord}
For each of these singular points, there is an analytic
 neighborhood isomorphic to
 a neighborhood of the origin of the set
\[
\{(X, Y, Z, t)\in\Bbb C^4\;|\; XY+tZ = 0\}\subset \Bbb C^3\times\Bbb C.
\]
\end{lem}
\proof
In fact, in terms of non-homogeneous coordinates
 $x = \frac{\bar x}{\bar w}, y = \frac{\bar y}{\bar w}$
 and $z = \frac{\bar z}{\bar w}$, the equation
 defining the degeneration becomes
\[
xyz+ tf = 0.
\]
We assume that the line $\ell_i$ is defined by $x = y = 0$, and 
 one of the singular point on it is given by $z = \alpha$, where $\alpha\neq 0$.
Let us write $\tilde z = z-\alpha$.
Then $f$ can be written in the form
\[
f = \tilde z f_1 + x f_2 + y f_3, 
\]
 where $f_i$ are polynomials and $f_1$ has nonzero constant term.
Moreover, $f_1$ depends only on $\tilde z$ and
 $f_2$ depends only on $x$ and $\tilde z$.
Write $\tilde y = y(\tilde z+\alpha) + tf_2$.
The equation becomes
\[
x\tilde y + t(\tilde z f_1+ \tilde y f_4+ tf_5) = 0,
\]
 where $f_4 = \frac{f_3}{\tilde z+\alpha}$, $f_5 = -\frac{f_2f_3}{\tilde z+\alpha}$.
Put $X = x+ tf_4$, $Y = \tilde y$ and $Z = \tilde zf_1+tf_5$. \qed
\begin{rem}
In fact, the isomorphism between a neighborhood of a point
 in $\mathcal S$ and a neighborhood of the origin
 of $\{(X, Y, Z, t)\in\Bbb C^4\;|\; XY+tZ = 0\}$
 can be obtained simply by using a coordinate system given by
 the functions $x, y$ and $\frac{f}{z}$.
However, we choose the above coordinate functions to simplify the
 later arguments.
\end{rem}
We consider the following situation.
Let
\[
\varphi_0\colon \Bbb P^1\to X_0\subset\Bbb P^3
\]
be an embedding of a projective line, whose image is contained in 
 an irreducible component of $X_0$.
The problem is the following:
\[
\begin{array}{c}
 \text{{\it Determine when $\varphi_0$ can be lifted to $X_t$, $t\neq 0$.}}
\end{array}
\]
As is pointed out several times (\cite{L}, see also \cite{NS} in the context of tropical 
 curves),  for an immersed stable curve
 $\psi_0\colon C\to X_0$, if the image is away from the singularities of $\mathfrak X$,
 it is necessary for $\psi_0$ to satisfy the \emph{pre-log} condition 
 (Definition 4.3, \cite{NS}) to solve the above problem.
In our situation, the line is embedded in a component of $X_0$, so we have the 
 following necessary condition.
\begin{lem}\label{lem:log}
Let $\varphi_0$ be as above.
Let $\Bbb P_i^2$ be the component of $X_0$ to which the line is mapped.
Then, to solve the problem above, it is necessary to satisfy the condition:
\[
(\ast)\;\; \varphi_0(\Bbb P^1)\cap L\subset \mathcal S. 
\]\qed\\
\end{lem}
Since the set $\mathcal S$ is contained in $\cup_{i=1}^6 \ell_i^{\circ}$,
 one sees that $\varphi_0$ is torically transverse as a map to $\Bbb P_i^2$, 
 in the sense of Definition 4.1 of \cite{NS}.
\begin{rem}\label{rem:genericity}
For the condition $(\ast)$ to be satisfied, a condition must be imposed on $f$, 
 the defining polynomial of the quartic surface.
Namely, the set $\mathcal S\cap\Bbb P^2_i$ should contain 3 points which are
 collinear.
This imposes 1 dimensional condition on $f$, and corresponds to the fact that
 a generic quartic
 K3 surface does not contain an embedded $(-2)$-curve.
In the rest of this section, we take generic $f$ among those which satisfy this condition. 
In particular, we still assume that the singular locus $\mathcal S$ does not intersect 
 the torus fixed point set of $\Bbb P^2_i$.
%
%
%
\end{rem}

The condition $(*)$ is necessary, but not sufficient.
For example, consider the family defined by the equation
\[
\bar x\bar y\bar z\bar w+t(\bar x^4-\bar z^4
 -2\bar z\bar w^3-\bar w^4+\bar y\bar w^3+\bar y^4) = 0
\]
 using homogeneous coordinates.
The set $\mathcal S$ contains the points
\[
(\bar x, \bar y, \bar z, \bar w) = (1, 0, 0, 1),\;\; (0, 0, -1, 1),\;\; (1, 0, 1, 0) 
\]
 so the line
\[
\bar x-\bar z-\bar w = 0,\;\; \bar y= 0 
\]
satisfies the condition $(*)$.
However, by direct calculation, one sees that there is no first order lift of this line.
The rest of this section is devoted to the study of the obstruction to the 
 existence of lifts of lines satisfying $(*)$.

\subsection{Calculation of the cohomology class}
As is usual with deformation theory, the possible
 obstructions are represented by suitable cohomology 
 classes.
However, 
 whether these obstructions really matter or not can not be seen just by 
 calculating
 the cohomology classes.
This is a non-linear problem, and solved by calculating {\it Kuranishi maps}.
We calculate the cohomology classes in this subsection, and 
 calculate Kuranishi maps in the next subsection.

\subsubsection{Identification of the log normal sheaf}\label{subsec:ns}
First we study the local structure around the singular points.
As noted above, the space $\mathfrak X$ has a neighborhood 
 $\mathfrak U$ which is analytically isomorphic to 
 a neighborhood of the origin $(X, Y, Z, t) = (0, 0, 0, 0)$ in 
 the variety defined by the equation $XY+tZ = 0$ at each singular point.
We consider an embedded line in $\Bbb P_i^2\subset X_0$
 which 
\begin{itemize}
 \item is torically transverse as a map to $\Bbb P_i^2$ and
 \item intersects the toric boundary at the singular points of $\mathfrak X$,
\end{itemize}
 as above.

From the proof of Lemma \ref{lem:coord}, 
 we see that locally the line is given by the
 equations
\[
aX+bZ +ZT(Z) = 0,\;\; Y = t = 0,
\]
 where $a, b$ are generic complex numbers.
In particular, neither $a$ nor $b$ is zero.
$T(Z)$ is a convergent series around $Z=0$ whose constant term is zero.
Thus, the line satisfies the following condition:
\[
(\diamondsuit)\hspace{.2in} \text{The line intersects the variety
 defined by $Z =Y = t = 0$ transversally on $\Bbb P^2_i$
 at $O$,}
\]
 here $O$ is the point $(X, Y, Z, t) = (0, 0, 0, 0)$.

We take a local affine coordinate $S$ on $\Bbb P^1$ so that 
\[
\varphi_0^*X = S
\]
 holds.
Here we write the composition 
\[
i_0\circ\varphi_0\colon \Bbb P^1\to X_0\to\mathfrak X
\]
 by the letter $\varphi_0$, for brevity.
Now we study how the local lifts of $\varphi_0$ are described using these coordinates.

Note that the variety $XY+tZ = 0$ has a natural structure of a toric variety
 over $Spec\;\Bbb C[t]$.
So it has a standard log structure which is log smooth over
 $Spec\;\Bbb C[t]$ ($Spec\;\Bbb C[t]$ is also equipped with the standard log structure
 as a toric variety).
Also, we put a log structure on $\Bbb P^1$ associated to the divisor 
 $z_0 = \{S = 0\}$.
The ghost sheaf of the log structure on $XY+tZ = 0$
 is isomorphic to the monoid
\[
\Bbb Z\langle\sigma_X, \sigma_Y, \sigma_Z, \sigma_t\rangle/
 (\sigma_X+\sigma_Y -\sigma_Z-\sigma_t).
\]
The ghost sheaf of the log structure on $\Bbb P^1$ in a neighborhood of
 a point in the inverse image of the set $\mathcal S$ is isomorphic to
\[
\Bbb Z\langle\tau_S, \tau_t\rangle,
\]

Then by the condition $(\diamondsuit)$,
 the map $\varphi_0$ can be equipped with a structure of 
 a map between log schemes with the following properties:
\begin{itemize}
\item The composition of $\varphi_0$ with the projection to
 $Spec\;\Bbb C[t]$ is log smooth.
\item The map $\varphi_0$
 is strict away from the inverse image of the singular
 set $\mathcal S$ of $\mathfrak X$.
\item The map $\varphi_0$ induces a map between ghost sheaves
 around the inverse image of $\mathcal S$ by
\[
\sigma_X\mapsto \tau_S,\;\; \sigma_Y\mapsto 0,\;\;
\sigma_Z\mapsto \tau_S,\;\; \sigma_t\mapsto \tau_t.
\]
\end{itemize}
The log structure 
 on $\varphi_0$ as a map between log schemes
 satisfying these conditions is uniquely determined up to isomorphisms.
\begin{rem}\label{rem:log}
It is important to note that we do not need a log structure 
 on whole $\mathfrak X$, but only on a neighborhood of the image of
 $\varphi_0$, to consider deformations of $\varphi_0$.
The log structures coming from the toric structure of $XY+tZ = 0$
 on neighborhoods of singular points of $\mathfrak X$ may not extend 
 to whole $\mathfrak X$.
However, under the condition $(\diamondsuit)$ 
 at each intersection between the line and the toric divisors,
 we can take the neighbourhood $\mathfrak U$ sufficiently small
 so that 
 the line $\varphi_0(C_0)$ does not intersect the variety
 $Z =Y = t = 0$ away from these intersections (note that 
 this variety is defined only on $\mathfrak U$).
Thus, on a suitable neighborhood of the image of
 $\varphi_0$, the log structures
 induced from the toric log structure of $XY+tZ = 0$
 around the singular locus of $\mathfrak X$
 extend by putting the pullback log structure
 from $Spec\;\Bbb C[t]$ away from the
 singular points.
\end{rem}

The log tangent sheaf $\Theta_{\mathfrak U}$
 of $\mathfrak U$ is locally free and generated by
 the sections
\[
X\partial_X+t\partial_t,\;\; Y\partial_Y+t\partial_t,\;\; Z\partial_Z-t\partial_t.
\]
There is a following relation between the log cotangent vectors:
\[
\frac{dX}{X}+\frac{dY}{Y}-\frac{dZ}{Z}-\frac{dt}{t} = 0,
\]
 which becomes
\[
\frac{dX}{X}+\frac{dY}{Y}-\frac{dZ}{Z}= 0
\]
 when restricted to $X_0$.
In particular, the subsheaf $\Theta_{X_0}$ restricted to 
 $X_0\cap \mathfrak U$ is generated by the sections
\[
X\partial X-Y\partial_Y, \,\, Y\partial_Y+Z\partial_Z.
\]
 

Since the map $\varphi_0$ factors through $X_0$, 
 the natural map from $\Theta_{\Bbb P^1}$ to $\varphi_0^*\Theta_{\mathfrak X}$
 is locally given by
\begin{equation}\label{eq:1}
S\partial_S\mapsto (X\partial_X-Y\partial_Y)
 +(1+\tilde T(Z))(Z\partial_Z+Y\partial_Y)
\end{equation}
 around $z_0$.
Here $\tilde T(Z)$ is a convergent series around $Z = 0$ whose constant
 term is zero.

Since $\mathfrak X\to \Bbb C$ 
 (restricted to a neighborhood of $\varphi_0(C_0)$, 
 see Remark \ref{rem:log})
 is log smooth, there is always a local lift of $\varphi_0$.
The obstruction to the existence of a
 lift of $\varphi_0$ is given by the first cohomology of the pull-back 
 $\varphi_0^*\Theta_{X_0}$ (\cite{KK}), 
 here $\Theta_{X_0}$ is the log tangent sheaf of $X_0$ with
 the log structure induced from that of $\mathfrak X$.
If a lift exists, the space of such lifts is a torsor of the zeroth cohomology
 of the same sheaf.
However, since the domain $\Bbb P^1$ (with 3 special points) has no moduli, 
 it suffices to replace $\varphi_0^*\Theta_{X_0}$ by the log
 normal sheaf for the calculation of the obstruction.
Now we calculate this log normal sheaf.

Let 
\[
j_0\colon \Bbb P^2_i\to \mathfrak X
\]
 be the inclusion.
In the neighborhood $\mathfrak U\cap \Bbb P_i^2$ of $\varphi_0(z_0)$, 
 we induce a log structure 
 from that of $\mathfrak U$ introduced above.
Let $\Theta_{\mathfrak U\cap\Bbb P_i^2}$ be the log tangent sheaf 
 of this log structure.

We remark that the coordinates $X, Z$ are not the standard coordinates 
 of an affine part of $\Bbb P^2_i$, which we wrote by $x, z$.
The log tangent sheaf $\Theta_{\Bbb P_i^2}^{st}$
 associated to the standard toric structure on 
 $\Bbb P^2_i$ is free of rank 2 and generated by
 $x\partial_x$ and $z\partial_z$.
The point $\varphi_0(z_0)$ has the parameter
 $(x, z) = (0, \alpha)$, $\alpha\neq 0$.
So on $\mathfrak U\cap\Bbb P_i^2$, the log tangent sheaf is generated by
 $x\partial_x$ and $\partial_z$.
Concerning the relation between the coordinates
 $X, Z$ and $x, z$, we note the following observation, which follows from 
 simple calculation.
\begin{lem}\label{lem:cood}
We have the following equalities between (log) tangent vectors
 on $\mathfrak U\cap\Bbb P_i^2$.
\begin{equation}\label{eq:2}
x\partial_x = (X\partial_X-Y\partial_Y) + \frac{x\partial_xf_1}{f_1}(
 Z\partial_Z+Y\partial_Y),
\end{equation}
\begin{equation}\label{eq:3}
\partial_{z}= \left(\frac{1}{\tilde z}+\frac{\partial_zf_1}{f_1}\right)
 (Z\partial_Z+Y\partial_Y),
\end{equation}
 where $\tilde z = z-\alpha$.\qed
\end{lem}
Note that the coefficient of 
 the right hand side of 
 (\ref{eq:3}) has the first order pole at $\varphi_0(z_0)$.
Thus, we see the following.
\begin{cor}
Using the coordinate $x$ and $z$, 
 the sheaf $\Theta_{\mathfrak U\cap \Bbb P_i^2}$ is generated by
 $x\partial_x$ and $\tilde z\partial_z$.\qed
\end{cor}
The map from $\Theta_{\Bbb P^1}$ to $\varphi_0^*\Theta_{\mathfrak X}$
 factors through the sheaf $\varphi_0^*\Theta_{\mathfrak U\cap \Bbb P_i^2}$.
Using the coordinate $x, \tilde z$, the line is given by the equation
\[
ax+b\tilde z = 0
\]
 as a subvariety of $\Bbb P^2_i$.
Near the point $z_0$, 
 the natural 
 map from $\Theta_{\Bbb P^1}$ to $\varphi_0^*\Theta_{\Bbb P_i^2}^{st}$
 is given by
\[
S\partial_S\mapsto x\partial_x+ \tilde z\partial_{z}.
\]
(here again we abused the notation by representing
 the map $\Bbb P^1\to \Bbb P_i^2\subset X_0$ by the same letter $\varphi_0$).
\begin{rem}
The image $x\partial_x+ \tilde z\partial_{z}$ is the same as
 the right hand side of (\ref{eq:1}).
From this and Lemma \ref{lem:cood}, we can represent $\tilde T(Z)$ in
 terms of the coordinates $x$ and $\tilde z$, though we do not use it.
\end{rem}

So in the neighborhood 
 $U_{z_0} = \varphi_0^{-1}(\mathfrak U\cap \Bbb P_i^2)$ of $z_0$, 
 the restriction $\mathcal N_{U_{z_0}}$ of the log normal sheaf
 $\varphi_0^*\Theta_{X_0}/\Theta_{\Bbb P^1}$ to $U_{z_0}$ is
 given by
\[
\mathcal N_{U_{z_0}} := \mathcal O_{U_{z_0}}
 \langle x\partial_x, \tilde z\partial_{z}\rangle/(x\partial_x+\tilde z\partial_{z}).
\]
 
On the other hand, the usual (non-log) normal sheaf of the 
 line $\{ax+b\tilde z = 0\}$ in $\Bbb P^2$ is,  in the neighborhood 
 $U_{z_0}$ of $z_0$, 
 given by the quotient
\[
\mathcal Q_{U_{z_0}} = \mathcal O_{U_{z_0}}\langle \partial_x, \partial_{z}\rangle/
 (b\partial_x- a\partial_{z}).
\] 
It is easy to see the following.
\begin{lem}
There is a natural map from $\mathcal N_{U_{z_0}}$
 to $\mathcal Q_{U_{z_0}}$
 which sends
 $x\partial_x$ to $x\partial_x$ and $\tilde z\partial_z$ to $\tilde z\partial_z$. 
\end{lem}
\proof
Note that since we have $ax+b\tilde z = 0$, 
\[
x\partial_x+\tilde z\partial_{z} = \frac{x}{b}(b\partial_x-a\partial_{z}).
\]
The assertion is clear from this. 
\qed\\
 
This is isomorphic except 
 the point $z_0 = \{S = 0\}$, and at this point,
 the quotient is the skyscraper
 sheaf $\Bbb C_{z_0}$.
It is clear that except at the singular point $z_0$, the log normal sheaf
 $\mathcal N_{\Bbb P_i^2/\Bbb P^1}$ is naturally isomorphic to the usual
 normal sheaf.
There are 3 singular points which are locally analytically isomorphic.
It follows from the natural isomorphism between 
 the sheaf $\mathcal N_{U_{z_0}}$
 and the usual 
 normal sheaf  $\mathcal Q_{U_{z_0}}$ except at singular points like $z_0$, 
 the 3 sheaves constructed as $\mathcal N_{U_{z_0}}$
 around the 3 singular points glue into a 
 global sheaf on $\Bbb P^1$
 which we write by $\mathcal N_{\Bbb P_i^2/\Bbb P^1}$.
This is the log normal sheaf $\varphi_0^*\Theta_{X_0}/\Theta_{\Bbb P^1}$
 we want.

The following is clear from definition.
Let $\mathcal Q$ be the usual normal sheaf of the map 
 $\varphi_0: \Bbb P^1\to \Bbb P_i^2$.
\begin{lem}\label{lem:N}
We have the exact sequence of sheaves
\[
0\to \mathcal N_{\Bbb P_i^2/\Bbb P^1}\to \mathcal Q\to \Bbb C_{p_1}\oplus
 \Bbb C_{p_2}\oplus\Bbb C_{p_3}\to 0,
\]
 here $p_i$ are the points on $\Bbb P^1$ mapped by $\varphi_0$ to 
 the set $\mathcal S$, the locus of singularities of $\mathcal X$.\qed
\end{lem}
The usual normal sheaf $\mathcal Q$ is isomorphic to 
 $\mathcal O_{\Bbb P^1}(1)$, so the log normal sheaf 
 $\mathcal N_{\Bbb P_i^2/\Bbb P^1}$ is isomorphic to $\mathcal O_{\Bbb P^1}(-2)$.

\subsubsection{Calculation of the (dual) obstruction class}
Since the sheaf $\mathcal N_{\Bbb P^2_i/\Bbb P^1}$ is isomorphic to
 $\mathcal O_{\Bbb P^1}(-2)$, 
 the cohomology groups of $\mathcal N_{\Bbb P^2_i/\Bbb P^1}$ are
\[
H^0(\mathcal N_{\Bbb P^2_i/\Bbb P^1}) = 0,\;\;
 H^1(\mathcal N_{\Bbb P^2_i/\Bbb P^1}) = \Bbb C.
\] 
That is, there is an obstruction class for deforming the map $\varphi_0$.
For the calculation of the Kuranishi map, we need to represent this class
 in an effective way.

By construction, $\mathcal N_{\Bbb P^2_i/\Bbb P^1}$ is the sheaf of sections of
 the usual normal sheaf $\mathcal Q$ which have single zeroes at
 3 points $p_1, p_2, p_3$ mapped to the singular locus $\mathcal S$.
By Serre duality, 
\[
H^1(\mathcal N_{\Bbb P^2_i/\Bbb P^1})
 \cong 
 (H^0((\mathcal N_{\Bbb P^2_i/\Bbb P^1})^{\vee}
  \otimes\omega_{\Bbb P^1}))^{\vee},
\]
 here $\omega_{\Bbb P^1}$ is the canonical sheaf of $\Bbb P^1$.
By the above description, 
 $(\mathcal N_{\Bbb P^2_i/\Bbb P^1})^{\vee}$ is the sheaf of 
 sections of the usual conormal sheaf $\mathcal Q^{\vee}$,
 where single poles are allowed at the 3 points $p_1, p_2, p_3$. 

Let us describe the space of global sections of the sheaves $\mathcal Q$ and
 $(\mathcal N_{\Bbb P^2_i/\Bbb P^1})^{\vee}$ geometrically.
The projective plane $\Bbb P_i^2$ has a natural structure of a toric variety.
Let $N = \Bbb Z^2$ be a lattice and
 $N_{\Bbb R}$ be the vector space where the fan describing $\Bbb P_i^2$ lies.
The fan has 3 rays generated by the integral vectors 
 $e_1 = (1, 0), e_2= (0, 1)$ and
 $(-1, -1)$ respectively.
Let $M_{\Bbb R}$ be the dual space of $N_{\Bbb R}$
 and $f_1, f_2$ be the dual basis of $e_1, e_2$.
First we observe the following.
\begin{lem}\label{lem:normal}
The space of global sections of the sheaf $\mathcal Q$ is naturally isomorphic to 
 $N_{\Bbb C} = N\otimes_{\Bbb Z}\Bbb C$.
\end{lem}
\proof
Recall that $\mathcal Q$ is the normal sheaf of the line $\varphi_0(\Bbb P^1)$
 in $\Bbb P_i^2$.
Let $\Theta_{\Bbb P^1}$ be the log tangent sheaf of $\Bbb P^1$ with
 respect to the log structure on $\Bbb P^1$ associated to the divisor
 $p_1+p_2+p_3$.
Then we have the following commutative diagram:
\[
\xymatrix{
0\ar[r] & \Theta_{\Bbb P^1}\ar[r]\ar[d]^{i_1} & \Theta^{st}_{\Bbb P^2_i}|_{\Bbb P^1}
 \ar[r]\ar[d]^{i_2} & \nu_{\Bbb P^2_i/\Bbb P^1}\ar[r]\ar[d]^{i_3} & 0\\
0\ar[r]& T_{\Bbb P^1}\ar[r]& T_{\Bbb P^2_i}|_{\Bbb P^1}\ar[r] & \mathcal Q \ar[r] & 0.}
\]
Here the rows are exact, and the sheaf $\nu_{\Bbb P^2_i/\Bbb P^1}$
 is defined by this sequence.
The sheaves $T_{\Bbb P^1}$ and $T_{\Bbb P_i^2}$ are the 
 usual (not log) tangent sheaves, and 
 the vertical maps $i_1, i_2$ are the natural injections.
The map $i_3$ is induced by $i_2$.
Then the sheaves $\nu_{\Bbb P^2_i/\Bbb P^1}$ and $\mathcal Q$ are both isomorphic to
 $\mathcal O_{\Bbb P^1}(1)$, and it is easy to see that the map $i_3$ is an injection.
Thus, $i_3$ is in fact an isomorphism.
So the space of global sections of $\nu_{\Bbb P^2_i/\Bbb P^1}$ and that of $\mathcal Q$
 are isomorphic, too.
Now by the upper row of the diagram, the spaces of global sections of 
 the sheaves $\Theta_{\Bbb P^2_i}|_{\Bbb P^1}$ and $\nu_{\Bbb P^2_i/\Bbb P^1}$
 are isomorphic, because the sheaf $\Theta_{\Bbb P^1}$ is isomorphic to
 $\mathcal O_{\Bbb P^1}(-1)$ and so all its cohomology groups vanish.
Since the space of the global sections of the sheaf  
 $\Theta_{\Bbb P^2_i}|_{\Bbb P^1}$ is naturally isomorphic to 
 $N_{\Bbb C}$, the lemma follows.\qed\\
 
Next, we prove the following. 
\begin{lem}\label{lem:conormal}
There is an isomorphism between
 the space of 
 global sections of the sheaf $(\mathcal N_{\Bbb P^2_i/\Bbb P^1})^{\vee}$
 and the 3 dimensional vector space whose bases are
 naturally identified with the vectors $f_1, f_2$ and $f_1-f_2$ in $M_{\Bbb R}$:
\[
H^0(\mathcal N_{\Bbb P^2_i/\Bbb P^1})^{\vee}\cong
 \Bbb C\langle f_1\rangle\oplus \Bbb C\langle f_2\rangle\oplus
  \Bbb C\langle f_1-f_2\rangle.
\]
\end{lem}
\proof
Recall that $(\mathcal N_{\Bbb P^2_i/\Bbb P^1})^{\vee}$ is the sheaf of 
 sections of the usual conormal sheaf $\mathcal Q^{\vee}$,
 where single poles are allowed at the 3 points $p_1, p_2, p_3$.
From this, if a section $s$ of $(\mathcal N_{\Bbb P^2_i/\Bbb P^1})^{\vee}$
 does not have a pole at any
 of the points $p_1, p_2, p_3$, then $s$ is 
 the zero section.
In particular, a section of $(\mathcal N_{\Bbb P^2_i/\Bbb P^1})^{\vee}$
 is determined by its value at the points $p_1, p_2, p_3$.

On the other hand, we saw in the proof of the last lemma that
 the usual normal sheaf $\mathcal Q$ is naturally isomorphic to
 the sheaf $\nu_{\Bbb P^2_i/\Bbb P^1}$,
 which is the quotient of $\Theta_{\Bbb P^2_i}|_{\Bbb P^1}$ by
 $\Theta_{\Bbb P^1}$.
Let the points $p_1, p_2, p_3$ be mapped to 
 the points 
\[
[0, \alpha_1, 1], [\alpha_2, 0, 1], [\alpha_3, \alpha_4, 0]
\]
 on $\Bbb P_i^2$, 
 respectively.
Since we take the polynomial $f$ defining the degenerating family of
 K3 surfaces generic, the complex
 numbers $\alpha_i$, $i = 1, 2, 3, 4$ are not zero. 

Then, under the inclusion $\Theta_{\Bbb P^1}\to
  \Theta_{\Bbb P^2_i}|_{\Bbb P^1}$
  and identifying the fibers of the sheaf $\Theta_{\Bbb P^2_i}|_{\Bbb P^1}$
  with $N_{\Bbb C} = \Bbb C\langle x\partial_x, z\partial_z\rangle$, 
 the fibers of $\Theta_{\Bbb P^1}$ over the points $p_1, p_2, p_3$ are
 mapped to the subspaces generated by
\[
x\partial_x, z\partial_z, x\partial_x+z\partial_z,
\] 
 respectively.
In terms of the basis $e_1, e_2$ of $N_{\Bbb C}$, these are generated by
\[
e_1, e_2, e_1+e_2,
\]
 respectively.
 
Thus, the fibers of the sheaf $\mathcal Q^{\vee}$ over the points $p_1, p_2, p_3$
 are naturally isomorphic to the subspaces
\[
f_2, f_1, f_1-f_2,
\]
 respectively.

Since a global section of $(\mathcal N_{\Bbb P^2_i/\Bbb P^1})^{\vee}$
 is determined by its behavior at the poles (that is, at the points $p_1, p_2, p_3$),
 the space of global section is isomorphic to the space
$\Bbb C\langle f_2\rangle\oplus \Bbb C\langle f_1\rangle\oplus
  \Bbb C\langle f_1-f_2\rangle$.\qed
\begin{rem}
Note that the isomorphism of Lemma \ref{lem:conormal}
 is not canonical, since the residues of a section of 
 $(\mathcal N_{\Bbb P^2_i/\Bbb P^1})^{\vee}$ at the poles are not 
 well defined.
On the other hand, the space of global sections of 
 the sheaf $(\mathcal N_{\Bbb P^2_i/\Bbb P^1})^{\vee}\otimes\omega_{\Bbb P^1}$
 is canonically isomorphic to a subspace of 
 $\Bbb C\langle f_2\rangle\oplus \Bbb C\langle f_1\rangle\oplus
  \Bbb C\langle f_1-f_2\rangle$, as we will see below.
\end{rem}
\begin{lem}
The space of 
 global sections of the sheaf $(\mathcal N_{\Bbb P^2_i/\Bbb P^1})^{\vee}
  \otimes\omega_{\Bbb P^1}$
  is naturally isomorphic to the 1 dimensional subspace of
  $\Bbb C\langle f_1\rangle\oplus \Bbb C\langle f_2\rangle\oplus
  \Bbb C\langle f_1-f_2\rangle$ as follows:
\[
H^0((\mathcal N_{\Bbb P^2_i/\Bbb P^1})^{\vee}
  \otimes\omega_{\Bbb P^1})\cong
  \{(a, b, c)\in \Bbb C^3\;|\;
   af_1+bf_2+c(f_1-f_2) = 0\}.
\]
\end{lem}
\proof
Note that $(\mathcal N_{\Bbb P^2_i/\Bbb P^1})^{\vee}
  \otimes\omega_{\Bbb P^1}$
 is the sheaf of sections of the conormal sheaf $\mathcal Q^{\vee}$
 valued 1-form where log poles are allowed at the points
 $p_1, p_2$ and $p_3$.
As we saw in the proof of Lemma \ref{lem:conormal}, 
 the fibers of the sheaf $\mathcal Q^{\vee}$ over the points $p_1, p_2$ and $p_3$ 
 are naturally identified with subspaces of the vector space $M_{\Bbb C}$
 spanned by the vectors $f_2, f_1$ and $f_1-f_2$, respectively.

Using this trivialization of the fibers of $\mathcal Q^{\vee}$, 
 residues of sections of the sheaf 
 $(\mathcal N_{\Bbb P^2_i/\Bbb P^1})^{\vee}
  \otimes\omega_{\Bbb P^1}$ at the points $p_1, p_2$ and $p_3$ are well-defined.
Thus, a section $\sigma$ of the sheaf
 $(\mathcal N_{\Bbb P^2_i/\Bbb P^1})^{\vee}
  \otimes\omega_{\Bbb P^1}$ is naturally represented by
  an element
\[
v = a\langle f_1\rangle+b\langle f_2\rangle +c\langle f_1-f_2\rangle
\]
 of the vector space
 $\Bbb C\langle f_1\rangle\oplus \Bbb C\langle f_2\rangle\oplus
  \Bbb C\langle f_1-f_2\rangle$.
Here $a\langle f_1\rangle, b\langle f_2\rangle$
 and $c\langle f_1-f_2\rangle$ 
 are the residues of the section at the points $p_2, p_1$ and $p_3$.

On the other hand, sections of $(\mathcal N_{\Bbb P^2_i/\Bbb P^1})^{\vee}
  \otimes\omega_{\Bbb P^1}$ make a natural pairing with sections of 
  the sheaf $\mathcal Q$ to give meromorphic 1-forms on $\Bbb P^1$
  which allow log poles at the points $p_1, p_2$ and $p_3$.
By Lemma \ref{lem:normal},
 the space of global sections of the sheaf $\mathcal Q$ is identified with 
 the vector space $N_{\Bbb C} = \Bbb C\langle e_1, e_2\rangle$.

Thus, the pairing of the sections is given by the natural pairing between 
 the vector spaces $N_{\Bbb C}$ and $M_{\Bbb C}$. 
So the pairing between 
 the section of $(\mathcal N_{\Bbb P^2_i/\Bbb P^1})^{\vee}
 \otimes\omega_{\Bbb P^1}$ represented by the vector
 $v = a\langle f_1\rangle+b\langle f_2\rangle +c\langle f_1-f_2\rangle$
 and the section of $\mathcal Q$ represented by the vector $e_1$
 is the meromorphic 1-form whose residues at the points $p_1, p_2, p_3$
 are $0, a, c$, respectively.  
Similarly, the pairing between 
 the section of $(\mathcal N_{\Bbb P^2_i/\Bbb P^1})^{\vee}
 \otimes\omega_{\Bbb P^1}$ represented by
 $v$ and 
 the section of $\mathcal Q$ represented by $e_2$
 is the meromorphic 1-form whose residues at the points $p_1, p_2, p_3$
 are $b, 0, -c$, respectively. 

Now by the residue theorem, 
 the sum of the residues of a meromorphic 1-form
 must be zero.
Thus,
\[
a+c = b-c = 0.
\]
This is equivalent to the condition that 
\[
af_1+bf_2+c(f_1-f_2) = 0.
\]
 in $M_{\Bbb C}$.

This subspace of $M_{\Bbb C}$ is 
 1 dimensional and the dimension of the space of
 global sections of $(\mathcal N_{\Bbb P^2_i/\Bbb P^1})^{\vee}
 \otimes\omega_{\Bbb P^1}$ is also 1 dimensional, so these spaces are isomorphic.
 \qed\\

Summarizing, we have the following.
\begin{thm}\label{thm:K3}
Let $\varphi_0\colon \Bbb P^1\to X_0\subset\mathfrak X$ be a line satisfying the condition
 {\rm ($\ast$)}
 of Lemma \ref{lem:log}.
The cohomology group representing the obstruction to the existence of a lift of 
 $\varphi_0$ is isomorphic to
 the 1 dimensional vector space $\Bbb C$, and 
 its dual space is canonically presented by the triplet of vectors
\[
v_1\in \Bbb C\cdot f_1,\;\; v_2\in\Bbb C\cdot f_2,\;\; v_3\in\Bbb C \cdot (f_1-f_2)
\]
 in $M_{\Bbb R}$ satisfying the condition
\[
v_1+v_2+v_3 = 0.
\]\qed
\end{thm}

\subsection{Remarks on Kuranishi maps}\label{subsec:Kuranishi}
It is now common to formulate abstract deformation theory in terms of
 differential graded Lie algebras (dgLa).
Deformation theory of a complex manifold or pair of them is a typical example.
The Kuranishi map, which is a nonlinear map between suitable
 vector spaces whose zero set is the moduli space of deformation, 
 can also be formulated in terms of dgLa (see for example \cite{Ma}).
However, such an abstract formulation contains
 Green's operator, 
 a solution of non-linear PDE, or some other object which plays the 
 role equivalent to the Green's operator,
 which is usually very difficult to write down
 in concrete problems.
 
Our suggestion in this paper is that degeneration technique is useful in 
 calculating the Kuranishi map.
For this purpose, we adopt \v{C}ech cohomology point of view.
Let us explain it in a simplified situation.
Let $X$ be a nonsingular complex variety and 
\[
i\colon Y\hookrightarrow X
\]
 be an inclusion
 of a nonsingular subvariety,
 and suppose we want to study the deformation of $i$ with $X$ and $Y$ fixed.
The first order deformation of $i$ is naturally identified with the section of 
 $H^0(i^*TX)$.
Taking suitable coordinate neighborhoods of $X$ and $Y$ in analytic category,
 the image of $i$ is locally presented by a vector valued function on $Y$.
A section of $H^0(i^*TX)$ gives a perturbation of these functions on each coordinate
 neighborhood, with (infinitesimal) parameter $t$.

Up to the first order of $t$, these perturbations on coordinate neighborhoods are
 compatible under the coordinate change
 (because we start from a section of the pull back of the
 tangent bundle), however for the second order and higher, 
 these local perturbations do not necessarily coincide with each other on
 the intersections of the coordinate neighborhoods.
These gaps, inductively with respect to the order of $t$, 
 compose the degree 1 cohomology classes.
In particular, we call the map from $H^0(i^*TX)$ to the cohomology 
 which is given by the coefficients of  
 order 2 in $t$,
\[
\kappa^{(1)}\colon H^0(i^*TX)\to H^1(i^*TX),
\]
 the \emph{Kuranishi map of order 1}.

For a section of $H^0(i^*TX)$ 
 for which the image of $\kappa^{(1)}$ is trivial
 (in other words, at a point in the zero locus of $\kappa^{(1)}$), we can 
 perturb the first order deformation of $Y$ (or the map $i$) in the second order of $t$
 on each coordinate neighborhood, so that
 they are now compatible up to second order of $t$.
The space of such a perturbation is a torsor over the group $H^0(i^*TX)$.

Let $i_1$ be a first order deformation of $i$ which allows a
 second order deformation.
Then the set of second order deformations
 of $i_1$ is (non canonically) parametrised by $H^0(i^*TX)$. 
Using these second order deformations, 
 we again calculate the first cohomology,
 which is now induced from the terms of order 3
 with respect to $t$.
Then we can define the Kuranishi map of order 2 for the map $i_1$:
\[
\kappa_{i_1}^{(2)}\colon H^0(i^*TX) \to H^1(i^*TX).
\]

Continuing this process, one has a sequence $\{i_j\}_{0\leq j\leq k}$
 of deformations of $i = i_0$ such that:
\begin{itemize}
\item The map $i_j$ is defined over $\Bbb C[t]/t^{j+1}$.
\item The map $i_{j+1}$ reduces to $i_{j}$ over $\Bbb C[t]/t^{j+1}$.
\end{itemize}
Then, to the map $i_k$, we can define the Kuranishi map of 
 order $k$:
\[
\kappa_{i_{k-1}}^{(k)}\colon H^0(i^*TX) \to H^1(i^*TX).
\]
If this map has zero in its image, then some $k$-th order deformation of
 $i_{k-1}$ (not necessarily $i_k$) has a $k+1$-th order deformation.
 
If we can continue this process infinitely many times, it defines
 a formal deformation of $i$.
Then in a favorable situation ($X$ and $Y$ are projective, for example)
 implicit function theorem \cite{A}
 assures the existence of an actual deformation.

In general situation, it is still difficult to perform these process.
However, thorough degeneration, it becomes possible to write down the
 local lifts rather explicitly, and it enables us to calculate the Kuranishi map.

We note that the setting we will study
 is a little different from the one explained above.
First, we consider a family of maps
\[
\mathfrak C\to \mathfrak X
\]
 over $Spf\Bbb C[[t]]$.
Here $\mathfrak C$ is a family of prestable curves and $\mathfrak X$ is 
 a toric degeneration of a variety $X$. 
Second, the domain of the Kuranishi map is not the 
 zeroth cohomology of the pull-back of the tangent sheaf, but 
 that of the appropriate normal sheaf (which we calculated in Subsection \ref{subsec:ns}). 
 
We calculate the Kuranishi map for lines in a $K3$ surface in the next subsection.
See \cite{N} for another example of the
 calculation of the Kuranishi map using degeneration.

\subsection{Calculation of the Kuranishi map}\label{subsec:b_1}
The Kuranishi map reveals whether 
 the class calculated in Theorem \ref{thm:K3}
 really obstructs the existence of 
 a lift of $\varphi_0$ or not.
We begin with studying the local model around the singularity.
Namely, consider the variety defined by the equation
\[
XY+tZ = 0
\]
 and a nonsingular curve $C$ in it given by a parametrization:
\[
X = s,\;\; Z = s\zeta(s),\;\; Y = t = 0.
\]
Here $\zeta$ is an analytic function and 
 $s$ is a coordinate on the curve $C$.
In our case, it suffices to take $C$ as a small open disk. 
Note that this curve intersects the singular point
\[
(X, Y, Z, t) = (0, 0, 0, 0).
\]

We study a general form of lifts of $C$ (over $\Bbb C[t]$).
They should have the form
\[
X = s+tf_1(s)+t^2 f_2(s)+\cdots,
\]
\[
Z = s\zeta(s)+tg_1(s)+t^2g_2(s)+\cdots,
\]
\[
Y = th_1(s)+t^2h_2(s)+\cdots,
\]
 where $f_i, g_i, h_i$ are analytic functions.

 

Taking $S = s+tf_1(s)+t^2f_2(s)+\cdots$ as a new coordinate on $C$
 (more precisely, take it as a coordinate on the product 
 $C\times Spec\;\Bbb C[t]/t^k$ for some $k$), 
 the above parametrization for $X, Y, Z$ becomes
\[
\begin{array}{ll}
X & = S,\\
Z & = (S-tf_1-t^2f_2-\cdots)\zeta(S-tf_1-t^2f_2-\cdots)
 +tg_1(S-tf_1-t^2f_2-\cdots)+\cdots\\
 & = S\zeta(S)+t\bar g_1(S)+t^2 \bar g_2(S)+\cdots,\\
Y &= t\bar h_1(S)+ t^2\bar h_2(S)+\cdots,
\end{array}
\]
 where $\bar g_i, \bar h_i$ are analytic functions.
By $XY + tZ = 0$, 
\[
\bar h_1(S) = -\zeta(S),\;\; S\bar h_{i+1}(S)+\bar g_i(S) = 0,\;\; i\geq 1.
\]
Thus, using an appropriate coordinate on $C$, we have the expression of 
 a general lift as
\[
X = S,
\]
\[
Z = -S\bar h_1(S)-tS\bar h_2(S)-t^2S\bar h_3(S)-\cdots,
\]
\[
Y = t\bar h_1(S)+ t^2\bar h_2(S)+\cdots.
\]
In particular, for any order of $t$, the lift intersects
 the subvariety defined by $X = Z = 0$ 
 at the parameter $S=0$.
This can also be seen from the discussion in the previous subsection,
 since the section of the normal sheaf must have a single zero at each singular 
 point, compared to the 'usual' normal sheaf of a line in $\Bbb P^2$. 

Now we return to our actual situation, where we consider a degeneration of
 K3 surfaces. 
Recall that in the proof of Lemma \ref{lem:coord}, we introduced a 
 change of coordinates
\[
X = x+tf_4,\;\; Y = \tilde y,\;\; Z = \tilde zf_1+tf_5
\]
 to locally bring the equation into the model case $XY+tZ = 0$.

In terms of the coordinates $x, y, \tilde z$, 
 the line in $X_0$ is represented as
\[
x = s,\;\; \tilde z = as,\;\; t = y = 0,
\]
 where $a$ is a nonzero complex number.
Lifts of this must have the form (after an appropriate affine change of the coordinates on 
 the domain curve)
\[
x = s,
\]
\[
\tilde z = as + t(a_1s+b_1)+t^2(a_2s+b_2)+\cdots,
\]
\[
y = t(c_1s+d_1)+t^2(c_2s+d_2)+\cdots.
\]

The relation between the parameters $s$ and $S$ of the line (near a singular point)
 is given by
\[
S = X = s+tf_4(s, t).
\]
Substituting this to $Z$, we have
\[\begin{array}{ll}
\tilde zf_1+tf_5
 & = (as+t(a_1s+b_1)+t^2(a_2s+b_2)+\cdots)f_1(s, t)+tf_5(s, t)\\
 & = -S\bar h_1(S)-tS\bar h_2(S)-t^2S\bar h_3(S)-\cdots\\
 & = -(s+tf_4(s, t))\bar h_1(s+tf_4(s, t))-t(s+tf_4(s, t))\bar h_2(s+tf_4(s, t))
  -\cdots.
\end{array}
\]
Note that $f_1(s, t), f_4(s, t)$ and $f_5(s, t)$ contains the variable $t$ because
 these are functions of $x, \tilde z, y$,  and $\tilde z$ and $y$ have the above form.

Comparing the coefficients of $s^i$, $i = 1, 2, \dots$, 
 we see
\[
asf_1(s, 0) = -s\bar h_1(s).
\]
So we have
\[
\bar h_1(s) = -af_1(s, 0).
\]
This determines $\bar h_1$.
From this, the coefficients $c_1$ and $d_1$ above are deterimned.

Now let us compare the coefficients of $t$.
We have
\[
b_1f_1(0)+f_5(0) = -f_4(0)\bar h_1(0).
\]
Since we know 
 $\bar h_1$, we can calculate $b_1$ from this (recall $f_1(0)\neq 0$).

Using this calculation, we can state the necessary and sufficient condition 
 for the existence of the first order lift.
Recall that the degeneration of K3 surface is defined by
 $\bar x\bar y\bar z\bar w+tf = 0$.
We use $\bar x$ etc. 
 instead of $x$ to distinguish the homogeneous (former)
 and inhomogeneous coordinates.
Let 
\[
\ell\colon \; \bar x+a\bar z+b\bar w = 0,\;\; \bar y = t = 0
\]
 be a line on $X_0$ satisfying the condition $(\ast)$.
Let $\Bbb P^2_{\bar y} 
 = \{t = \bar y = 0\}$ be the component
  of $X_0$ on which the line $\ell$ lies.
The intersection of $\ell$ and the set of the toric divisors
 of $\Bbb P_{\bar y}^2$
  is the 3 point set
 composed of
\[
\left(\frac{\bar x}{\bar w}, 
  \frac{\bar z}{\bar w}\right) =  \left(0, -\frac{b}{a}\right),\;\; \left(-b, 0\right),
\]
 and
\[
\left(\frac{\bar z}{\bar x}, \frac{\bar w}{\bar x}\right)
 = \left(-\frac{1}{a}, 0\right).
\]
Let us write these points by $p_1, p_2, p_3$.
At each of these intersections, we compute the term $b_1$ as above.
Let $b_{1, \bar x}, b_{1, \bar z}, b_{1, \bar w}$ be the value of $b_1$ at
 $p_1, p_2, p_3$ respectively.
\begin{thm}\label{thm:Kuranishi}
The line $\ell$ has the first order lift if and only if 
\[
(\star)\colon ab_{1, \bar x}-abb_{1, \bar w}-b_{1, \bar z} = 0
\]
 holds.
\end{thm}
\proof
As we saw in Subsection \ref{subsec:ns}, 
 the space of local first order lifts around $p_i$ is
 a torsor over
 the vector space of local sections of 
 $\mathcal N_{\Bbb P^2_{\bar y}/\Bbb P^1}$.
We write this vector space by $\mathcal V_{p_i}$.
Let $x = \frac{\bar x}{\bar w}$, 
 $y = \frac{\bar y}{\bar w}$ and $z = \frac{\bar z}{\bar w}$ as before.
The functions $x, y, z, t$ compose a local coordinate
 of the ambient space $\Bbb P^3\times \Bbb C$ in which the 
 total space $\mathfrak X$ of the degeneration lies.  
By the calculation above, at $p_1$, this torsor can be written in the form
\[
(\ast) \hspace{.2in} \partial_t+(c_1s+d_1)\partial_y+ 
 (a_1s+b_{1, \bar x})\partial_{z} + \mathcal V_{p_1}.
\]
 using these coordinates.
\begin{rem} 
This is just a formal expression, since the sum of the tangent vectors
 on $\Bbb P^3\times \Bbb C$ and the vectors in $\mathcal V_{p_1}$
 does not make sense.
However, this expression is still reasonable in the following sense.
Namely, regarding the curve $\ell$ as a map to $\Bbb P^3\times\Bbb C$, 
 its first order deformations canonically correspond
 to sections of the pull-back of
 the tangent sheaf of $\Bbb P^3\times\Bbb C$.
In particular, the lift 
 parametrized by
 $x = s, \tilde z = as+t(a_1s+b_{1, \bar x}), y = t(c_1s+d_1)$
 considered above corresponds to the section 
 $\partial_t+(c_1s+d_1)\partial_y+ 
 (a_1s+b_{1, \bar x})\partial_{z}$.
Other lifts which are contained in $\mathfrak X$
 (modulo automorphisms)
 are given by perturbing this, by elements of $\mathcal V_{p_1}$.
\end{rem} 
 
The same holds at $p_2$ and $p_3$.
The line $\ell$ has the first order lift if and only if the \v{C}ech 1-cocycle
 $\sigma\in H^1(\Bbb P^1, \mathcal N_{\Bbb P_{\bar y}^2/\Bbb P^1})$
 defined by the difference of the local lifts is cohomologically trivial.
In other words, the line $\ell$ has 
 the first order lift if and only if the coupling between this cohomology 
 class and the generator of the class of Theorem \ref{thm:K3}
 is 0.
 
Now consider the class calculated in Theorem \ref{thm:K3}.
One sees that the fibers of the
 invertible sheaf 
 $(\mathcal N_{\Bbb P_{\bar y}^2/\Bbb P^1})^{\vee}$ are spanned by
\[
\frac{d\bar z}{\bar z}-\frac{d\bar w}{\bar w},\;\;
\frac{d\bar w}{\bar w}-\frac{d\bar x}{\bar x},\;\;
\frac{d\bar x}{\bar x}-\frac{d\bar z}{\bar z},
\]
 at $p_1, p_2, p_3$ respectively.
These covectors can be naturally identified with the vectors
 $v_1, v_2, v_3$ of the vector space $M_{\Bbb R}$
 as in Theorem \ref{thm:K3}, and the generator of
 $H^0((\mathcal N_{\Bbb P_{\bar y}^2/\Bbb P^1})^{\vee}
 \otimes \omega_{\Bbb P^1})\cong\Bbb C$
 is uniquely determined by them.
Namely, the values of the section 
 which generates $H^0((\mathcal N_{\Bbb P_{\bar y}^2/\Bbb P^1})^{\vee}
 \otimes \omega_{\Bbb P^1})\cong\Bbb C$
 at these points are
\[
\left(\frac{d\bar z}{\bar z}-\frac{d\bar w}{\bar w}\right)\otimes
 \frac{ds_1}{s_1},\;\;
\left(\frac{d\bar w}{\bar w}-\frac{d\bar x}{\bar x}\right)\otimes
 \frac{ds_2}{s_2},\;\;
\left(\frac{d\bar x}{\bar x}-\frac{d\bar z}{\bar z}\right)\otimes 
 \frac{ds_3}{s_3},\;\;
\]
 where $s_i$ is a local coordinate of $\Bbb P^1$ at $p_i$, as in the calculation above.

This section couples with the local lifts given by $(\ast)$.
Using inhomogenous coordinates, we have 
\[
v_1 = \frac{d\bar z}{\bar z}-\frac{d\bar w}{\bar w}
 = \frac{d\frac{\bar z}{\bar w}}{\frac{\bar z}{\bar w}}
 = \frac{dz}{z}
 = -\frac{a}{b}dz.
\]
 at $p_1$.
Thus, around the point $p_1$, 
 the coupling of an element of $(\ast)$ and 
 $\left(\frac{d\bar z}{\bar z}-\frac{d\bar w}{\bar w}\right)\otimes
 \frac{ds_1}{s_1}$
 is given by $-\frac{a}{b}(a_1s_1+b_{1, \bar x})\otimes\frac{ds_1}{s_1}$
 plus the terms coming from the coupling of $-\frac{a}{b}dz$
 and the sections of $\mathcal V_{p_1}$.
Note that the residue does not depend on the
 latter part, and it is given by
\[
-\frac{a}{b}b_{1, \bar x}
\]

Similarly, from the points $p_2, p_3$, we obtain the residues
\[
\frac{1}{b}b_{1, \bar z},\;\; ab_{1, \bar w},
\] 
 respectively.
This computation gives the dual pairing between 
 $H^1(\Bbb P^1, \mathcal N_{\Bbb P_{\bar y}^2/\Bbb P^1})$
 and $H^0((\mathcal N_{\Bbb P_{\bar y}^2/\Bbb P^1})^{\vee}
 \otimes \omega_{\Bbb P^1})$ 
 (see \cite[Section 7]{F}).
Thus, the equation $(\star)$ is equivalent to the 
 vanishing of the sum of these 
 residues.
In other words, the equation $(\star)$ is equivalent to the 
 condition that the obstruction cohomology class
 $\sigma\in H^1(\Bbb P^1, \mathcal N_{\Bbb P_{\bar y}^2/\Bbb P^1})$
 vanishes.
This proves the proposition.\qed

\begin{rem}
Since the zeroth cohomology group
 of the log normal sheaf $\mathcal N_{\Bbb P^2_i/\Bbb P^1}$
 is zero and the first cohomology group is 1 dimensional, 
 the Kuranishi map 
 can be represented as
 a map from a point to $\Bbb C$, that is, it is given by the scalar
 $ab_{1, \bar x}-abb_{1, \bar w}-b_{1, \bar z}$. 
\end{rem}
This is the first order Kuranishi map.
So far, the coefficients $b_1, c_1$ and $d_1$ are determined.
When the Kuranishi map is zero,
 then we can determine $a_1$ from the values of $b_1$ at 2 singular
 points,
 and these coefficients determine the first order lift.
Then doing the same procedure for $b_2$,
 we can define the second order Kuranishi
 map, and when this is zero, we can find the second order lift, and so on.

\subsection{Example}
As an example of our calculation, we consider the degeneration given by
 the equation of the following form:
\[
\bar x\bar y\bar z\bar w+t(\bar x^4+\bar y^4-\bar z^4-\bar w^4
 +\bar y\bar w^3-2\bar z\bar w^3+a\bar x^2\bar w^2+b\bar x\bar z
 \bar w^2+c\bar x\bar w^3) = 0.
\]
Here $a, b, c, d$ are unfixed scalars.
This equation is chosen in the following way:
Consider a simpler equation
\[
\bar x\bar y\bar z\bar w+t(\bar x^4+\bar y^4-\bar z^4-\bar w^4
 +\bar y\bar w^3-2\bar z\bar w^3) = 0.
\]
This is the equation we gave after Remark \ref{rem:genericity}.
The singular locus of this degeneration contains 3 points
\[
(\bar x, \bar y, \bar z, \bar w) = (1, 0, 0, 1),\;\; (0, 0, -1, 1),\;\; (1, 0, 1, 0), 
\]
 which lie on the line
\[
\bar x-\bar z-\bar w = 0,\;\; \bar y = t = 0.
\]
However, it is easy to see (by calculating the Kuranishi map, or by calculating directly)
 that this line does not have the first order lift.
So we add terms 
 $a\bar x^2\bar w^2+b\bar x\bar z\bar w^2
 +c\bar x\bar w^3$, and study when 
 the above line has the first order lift.

First of all, the singular locus of the total space $\mathfrak X$ must contain
 the above 3 points.
This implies the condition
\[
a+c = 0.
\]
So the equation becomes
\[
\bar x\bar y\bar z\bar w+t(\bar x^4+\bar y^4-\bar z^4-\bar w^4
 +\bar y\bar w^3-2\bar z\bar w^3+a\bar x^2\bar w^2
 +b\bar x\bar z\bar w^2-a\bar x\bar w^3) = 0.
\]
Now we assume $w\neq 0$ and inhomogenize the equation:
\[
xyz+t(x^4+y^4-z^4-1+y-2z+ax^2+bxz-ax) = 0.
\]
Consider the singular point $(x, y, z) = (1, 0, 0)$.
The parametrization is 
\[
x = s+1,\;\; z = s,\;\; y = 0.
\]
The functions $f_i$ are given as follows:
\[
f_1 = x^3+x^2+(1+a)x+1,
\]
\[
f_2 = -z^3-2+bx,
\]
\[
f_3 = 1+y^3.
\]
Then 
\[
f_4 = \frac{1+y^3}{x},\;\; f_5 = \frac{(1+y^3)(z^3+2-bx)}{x}. 
\]
Since
\[
b_1 = \frac{f_1(0)f_4(0)-f_5(0)}{f_1(0)}, 
\]
 we have
\[
b_1 = \frac{4+a-2+b}{4+a} = \frac{a+b+2}{4+a}.
\]
The case when $4+a = 0$ corresponds to the situation where 
 two of the points in the singular set $\mathcal S$
 of $\mathfrak X$ (see Section \ref{sec:K3}) merges.
So under our assumption that $f$ is generic, this does not happen.

Similarly, at $(x, y, z) = (0, 0, -1)$, one calculates
\[
b_1 = \frac{a+b-2}{2}.
\]
At the remaining singularity which lies at $x, z\to \infty$, one sees $b_1 = 0$.
Thus, the Kuranishi map is given by 
\[
 \frac{a+b+2}{4+a}+\frac{a+b-2}{2}.
\]
The condition for the vanishing of the Kuranishi map is
\[
a^2+ab+4a+6b-4 = 0.
\]

On the other hand, we can directly calculate the liftability
 condition by substituting 
\[
x = s,
\]
\[
\tilde z = as + t(a_1s+b_1)+t^2(a_2s+b_2)+\cdots,
\]
\[
y = t(c_1s+d_1)+t^2(c_2s+d_2)+\cdots.
\]
 into the equation.
Then from the coefficients of $ts^i$, we have
\[
c_1 = -4,\;\; d_1 = 2-a-b.
\]
Using this, we see from the coefficients of $t^2s^i$, 
\[
a_1 = 0,\;\; 2b_1+2-a-b = 0,\;\; (a+6)b_1+4 = 0.
\]
Solving this, we again have
\[
a^2+ab+4a+6b-4 = 0.
\]
\begin{rem}
In general, calculating the coefficients $a_i, b_i, c_i, d_i$
 for larger $i$ is quite hard, 
 though the method in this section reduces the amount of the calculation.
However, knowing that we can calculate 
 the Kuranishi map in this way is important.
In some situations, it happens that 2
 quite different curves have the same
 description of the Kuranishi maps, and one of these curves has vanishing
 obstruction by an obvious reason.
In this case, the obstruction of the other curve must vanish, which gives
 a nontrivial deformability result.
In \cite{N5}, we show the existence of infinitely many rational curves
 on quartic K3 surfaces which belong to an open dense subset of the 
 moduli space based on this method. 
\end{rem}

\begin{rem}
Each component of $X_0$ is $\Bbb P^2$ with 4
 distinguished points (singular locus
 of $\mathfrak X$) on each of the toric divisors.
We considered a line through 3
 of these points.
There is a unique cubic
 curve through the remaining 9 points on $\Bbb P^2$.
If the line is liftable, then the 
 intersection of a hyperplane containing the lifted line
 with the K3 surface contains a cubic curve
 which lifts the given cubic curve on $X_0$,
 so the above calculation of the 
 Kuranishi map also calculates the Kuranishi map for the cubic curve.
\end{rem}

\subsection{Cubic surface}
\begin{figure}[h]
\includegraphics{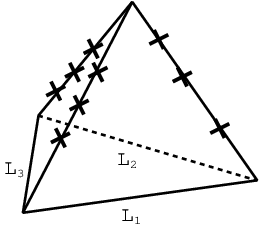}
\caption{}\label{fig:1}
\end{figure}

As an easy application of the calculation of the proof of Theorem \ref{thm:K3},
 we consider the 
 case of a cubic surface.
We consider a degeneration similar to the K3 case:
\[
xyz + tf = 0,
\]
 where $x, y, z$ are three of the standard 
 homogeneous coordinates of $\Bbb P^3$, 
 $t\in\Bbb C$ is the parameter for the degeneration, and $f$ is a generic 
 homogeneous cubic polynomial.
Let $\mathfrak X\subset \Bbb P^3\times\Bbb C$ be the total space of the degeneration,
 and $X_0\subset\mathfrak X$ be the central fiber as before.
In this case, $X_0$ is a union of 3 projective planes glued along 
 3 projective lines:
\[
X_0 = \Bbb P_1^2\cup \Bbb P_2^2\cup \Bbb P_3^2,
\]
\[
 \Bbb P_i^2\cap \Bbb P_{i+1}^2
 = \ell_i, \;\; (i\in \Bbb Z/3\Bbb Z). 
\]
On each of these lines, there are 3 singular points
 of $\mathfrak X$.
See Figure \ref{fig:1}.

Note that the bottom is open.
Each singular point has the same local structure as in the K3 case. 
 
We consider the problem of which line in $X_0$ lifts to a general fiber of 
 $\mathfrak X$.
So let 
\[
\varphi_0\colon \Bbb P^1\to \Bbb P_i^2
\]
 be a line.
As in the K3 case, if $\varphi_0$ is liftable, then it is necessary that the intersection
 with $\ell_i$ must be at the singular points.
However, in this case, there is no restriction to the intersection with the bottom lines
 ($L_1, L_2, L_3$ in the figure).

The relevant normal sheaf $\mathcal N_{\Bbb P_i^2/\Bbb P^1}$ is
\[
\mathcal O_{\Bbb P^1}(1-2) = \mathcal O_{\Bbb P^1}(-1)
\]
 by the calculation in Subsection \ref{subsec:ns}.
In this case, the cohomology groups are 
\[
H^0(\mathcal N_{\Bbb P^2_i/\Bbb P^1}) = 0,\;\;
 H^1(\mathcal N_{\Bbb P^2_i/\Bbb P^1}) = 0.
\] 
In particular, the line is unobstructed.
Thus, any line which intersects $\ell_i$ at the singular points lifts, 
 and on each $\Bbb P^2_i$, there are $3\times 3 = 9$ such lines.
So there are $9\times 3 = 27$ liftable lines in total, which proves the well-known 
 result for the cubic surface.

\section{Counting lines in quintic Calabi-Yau hypersurfaces}\label{sec:CYline}
\subsection{Combinatorics of lines in the central fiber}\label{subsec:combi}
Here we review an idea to combinatorially count
 $(-1)$-curves in quintic Calabi-Yau hypersurfaces via degeneration,
 due to Sheldon Katz \cite{Katz0, Katz}.

Consider a degeneration of a generic quintic Calabi-Yau hypersurface
 defined by the equation
\[
z_0z_1z_2z_3 z_4+ tf = 0,
\]
 where $z_0, \dots, z_4$ are homogeneous coordinates of $\Bbb P^4$, $t\in\Bbb C$ is
  the parameter of the 
  degeneration and $f$ is a generic homogeneous quintic polynomial
 of $z_0, \dots, z_4$.
Let $\mathfrak X$ be the total space
 and 
\[
X_0 = \cup_{i = 1}^5 \Bbb P_i^3
\]
 be the central fiber.
Each $\Bbb P^3_i$ has a natural structure of a toric variety.
Let us consider $\Bbb P^3_1$.
In this case, the intersection of $\Bbb P^3_1$ with
 the singular locus 
 $\mathcal S$ of $\mathfrak X$ is given by the union of 4 quintic curves, 
 one for each toric divisor:
\[
\Bbb P^3_1\cap 
 \mathcal S = \cup_{j=2}^5C_j,\;\; C_j\subset \Bbb P_1^3\cap\Bbb P_{j}^3
\]
Moreover, for $j_1\neq j_2$, $C_{j_1}\cap C_{j_2}$ is the set of 5 points. 
Since $f$ is generic, we can assume that $C_j$ does not intersect the 
 torus fixed points of $\Bbb P_1^3$.

We consider the lifting problem for lines in $\Bbb P^3_1$.
As in the K3 case, 
 it is easy to see that
 a liftable line must intersect the toric divisor of $\Bbb P^3_1$ 
 only at $\Bbb P^3_1\cap \mathcal S$.
Since a line intersects each toric divisor only once
 and $C_j$s lie on the toric boundary, this condition is the same as the condition that
 the line intersects each $C_j$.

It is well-known that the number of lines in $\Bbb P^3$ which intersect 
 generic 4 lines is 2.
This corresponds to the fact that
 the degree of $Gr(2, 4)$ under the Pl\"ucker embedding is 2.
Since we are considering the number of lines
 which intersect 4 curves of degree 5, there are at least
\[
2\cdot 5^4
\]
 lines satisfying this condition, for generic $f$.
\begin{defn}
Let $\mathcal L$ be this set of lines in $\Bbb P_1^3$.
\end{defn}
\begin{rem}
Even if $f$ is generic, the configuration of the 4 quintic curves in $\Bbb P_1^3$
 (that is, $\Bbb P_1^3\cap \mathcal S$)
 is not generic in the space of the configurations
 of such curves.
So the number of lines might not be $2\cdot 5^4$, but there are at least
 $2\cdot 5^4$ lines, when counted with multiplicity.
\end{rem} 

These lines can be classified into 2 classes.
Let 
\[
\ell_{j_1j_2} = \Bbb P_1^3\cap \Bbb P_{j_1}^3\cap \Bbb P_{j_2}^3,\;\;
 j_1\neq j_2,\;\; j_1\neq 1, \;\; j_2\neq 1
\]
 be the toric subvariety of codimension 2 in $\Bbb P^3_1$.
Then the 2 cases are:
\begin{enumerate}
\item The lines which do not intersect $\ell_{j_1j_2}$ for any $j_1, j_2$. 
\item The lines which intersect some $\ell_{j_1j_2}$.
\end{enumerate}
Now we count the number of lines in class (2).
The lines in class (2) can be further divided into 2 subclasses:
\begin{enumerate}
\item[(2)-I] The lines which intersect 2 different $\ell_{j_1j_2}$s.
\item[(2)-II] The lines which intersect only 1 $\ell_{j_1j_2}$.
\end{enumerate}
Let  $p$ be a point in $\ell_{j_1j_2}\cap \mathcal S$.
Then the number of lines containing $p$ in the set $\mathcal L$
 is 25, which is the number of the intersection points of the quintic curves
 $C_{j_3}, C_{j_4}$, projected to $\Bbb P^2$ from $p$
 (it is easy to see that we can take $f$ so that the projected 
 quintic curves intersect transversely), see also the proof of Proposition
 \ref{prop:transv}.

There are 30 points like $p$, so this counts
\[
25\cdot 30
\]
 lines of class (2).
However, this doubly counts the lines of class (2)-I.
It is easy to see that the number of lines in class (2)-I is
\[
5^2\cdot 3, 
\]
 where 3 is the number of pairs of $\ell_{j_1j_2}$ which does not intersect.
Thus, there are
\[
25\cdot 27 = 675
\]
 lines of class (2).

We see that the number of the lines of class (1) is at least,
\[
2\cdot 5^4 - 25\cdot 27 = 575.
\]
This is the number we want, since 
\[
575\cdot 5 = 2875
\]
 is the known number of lines
 in a generic quintic Calabi-Yau hypersurface.

We prove the following, which justifies Katz's idea
 sketched in \cite{Katz}. 
\begin{thm}\label{thm:CY3}
For generic $f$, there are exactly 2875 distinct lines of class (1) and 
 each of them lifts in a unique way.
The lines of class (2) do not lift.
\end{thm}

\subsection{Proof of the theorem}
\subsubsection{Finiteness and transversality}
First we prove that the number of lines in $\Bbb P^3_1$ which satisfy the 
 incidence conditions is finite.
This is immediate when the incidence conditions (4 quintic curves) are generic
 and the corresponding Schubert calculus gives the expected answer.
However, since our incidence conditions are in very special position, we
 prove the finiteness by low-tech argument.
\begin{lem}\label{lem:finite0}
The number of lines in $\Bbb P_1^3$ whose intersection with the 
 toric boundary is contained in the set $\mathcal S$ is finite for any generic
 quintic polynomial $f$.
\end{lem}
\proof
As we saw in the previous subsection,
 for general $f$, the number of lines of class (2) is finite.
A line of class (1) can be described as follows.
Recall that
 we write the intersection of $\Bbb P^3_1$ and the singular locus
  $\mathcal S$
 of $\mathfrak X$ by
\[
\Bbb P^3_1\cap \mathcal S = \cup_{j=2}^5 C_j.
\]
Fix a point $p$ in $C_2$ which is contained in the open orbit of the
 toric boundary of $\Bbb P^3_1$ 
 and consider the projection
\[
\pi_p\colon \Bbb P^3_1\to \Bbb P^2
\]
 from $p$.
The images of $C_3, C_4, C_5$ give
 3 quintic curves in $\Bbb P^2$.
A line through $p$ satisfying the incidence conditions corresponds to a
 point in the intersection
\[
C_3\cap C_4\cap C_5.
\]
So the number of lines satisfying the incidence
 conditions which moreover intersect $p\in C_2$
 is finite and bounded by 25,
 except when the images of the curves
 coincide.
But this does not happen for any $p$ when $f$ is generic.

On the other hand, the set of those $p\in C_2$ for which there is a 
 line through it satisfying the incidence conditions
 is a closed analytic subset, which is not the whole $C_2$ for generic $f$.
So the cardinality of the set of such $p$ is also finite.
This proves the lemma.\qed

\begin{prop}\label{prop:transv}
For a general quintic polynomial $f$, there
 are exactly 575 lines of class (1), and 675 lines of class (2)
 in each $\Bbb P_i^3$.
\end{prop}
\proof
First, we note that the moduli space
 of plane quintic curves is 20 dimensional, so 
 even after fixing 15 points, namely 
 the position of the intersections with 
 the toric divisors of $\Bbb P^2$, we still have freedom to perturb
 the curve (these points 
 of intersection cannot be taken freely. Only 14 points among them can be 
 moved freely, but it does not matter to later argument).
 
Let $H_i$, $i = 2, 3, 4, 5$ be the hyperplanes in $Gr(2, 4)$
 composed of the lines in $\Bbb P^3_1$ intersecting the curves $C_i$, 
 $i = 2, 3, 4, 5$, respectively.
Since the transversality is an open condition, and by Lemma 
 \ref{lem:finite0} the number of lines satisfying the incidence conditions
 is finite, we only need to prove that when we pick up a particular
 line satisfying the incidence conditions, the intersection of 
 $H_i$ at the point $x$ 
 in $Gr(2, 4)$ corresponding to this line is transverse.
This transversality means the following condition.
Namely, let $v_i$ be a vector in the cotangent space $T^*_xGr(2, 4)$
 which annihilates the subspace $T_xH_i$.
Then the vectors $v_i$, $i = 2, 3, 4, 5$ generates $T^*_xGr(2, 4)$.

First we study the case of a line of class (1).
Note that since a line $\ell$ in
 $\Bbb P^3_1$ is determined by two points $a, b$ on it, 
 the tangent space of it as a point $x$ 
 in $Gr(2, 4)$ can be identified with the
 direct sum of the normal spaces of $\ell$ at the points $a$ and $b$.
Take $a, b$ as the intersections between $\ell$ and the curves
 $C_2, C_3$, and let $N_a\ell$, $N_b\ell$ be the normal space of 
 $\ell$ at these points.
  
Take a basis $\{w_{2, 1}, w_{2, 2}\}$ of the space $N_a\ell$ so that
 $w_{2, 1}$ is tangent to $C_2$.
Similarly, take a basis $\{w_{3, 1}, w_{3, 2}\}$ of the space
 $N_b\ell$ so that $w_{3, 1}$ is tangent to $C_3$.
As we remarked above, the vector space whose basis is 
 $\{w_{2, 1}, w_{2, 2}, w_{3, 1}, w_{3, 2}\}$ can be identified with
 $T_xGr(2, 4)$.

Now consider the matrix whose components are
 given by the values of the pairing between 
 the vectors $\{w_{i, j}\}$ and $\{v_k\}$.
It is easy to see that it has the form
\[
\begin{pmatrix}
0 & 0 & r & v\\
p & 0 & s & w\\
0 & 0 & t & y\\
0 & q & u & z
\end{pmatrix},
\]
 where its $k$-th column is given by the pairing between
 $\{w_{i, j}\}$ and $\{v_{k+1}\}$.
Here 
 $p$ and $q$ are nonzero.
Moreover, it is easy to see that for generic $f$, the matrix
\[
\begin{pmatrix}
r & v\\
t & y
\end{pmatrix}
\]
 is regular.
Thus, the vectors $v_2, v_3, v_4, v_5$ span the cotangent space 
 $T^*_xGr(2, 4)$.

The calculation for the lines of class (2) is similar, and we omit it.
So for generic $f$, 
 the 4 hypersurfaces $H_i$ in $Gr(4, 2)$, intersects transversally, 
 and gives expected $2\cdot 5^4$ intersections.
So the numbers of lines of classes (1) and (2) in 
 $\Bbb P^3_1$ are also as expected, 
 namely, 575 and 675.\qed

\subsubsection{Lines intersecting lower dimensional toric strata.}\label{subsec:(2)}
In this section, we consider the lifting of the lines of class (2) in Subsection
 \ref{subsec:combi},
 and see that these lines do not lift.
 
We consider a line of class (2)-II.
In this case, the line is contained in the closure of
 an orbit of a 2 dimensional subtorus of the
 big torus acting on $\Bbb P^3_1$.
Let us assume that the line intersects 
 the stratum $\ell_{12}$ (see Subsection 
 \ref{subsec:combi} for the notations).

As before, we consider the local model of the degeneration.
In this case, 
 mimicking the calculation in 
 the proof of Lemma \ref{lem:coord}, 
 one sees the following.
Namely, locally around the intersection of the line and
 $\ell_{12}$ (which should be a singular point in the total space of the 
 degeneration), the total space $\mathfrak X$ is 
 analytically isomorphic to the analytic subset of 
 $\Bbb C^4\times \Bbb C$
 defined by the equation
\[
xyz+ tw = 0,
\]
 where $\{x, y, z, w\}$ is the coordinate of $\Bbb C^4$ and 
 $t$ is the parameter of the degeneration.
We assume that the component $\Bbb P^3_1$ corresponds to 
\[
z = 0,\;\; t= 0.
\]
In this expression, the stratum $\ell_{12}$ corresponds to  
\[
x=y=z=t=0
\]
 and the point 
\[
x=y=z=w=y=0
\]
 corresponds to the unique singular point in a suitable open subset of
  $\mathfrak X$ lying on $\ell_{12}$.
Let  us write this point by $p$.
Our line should intersect this point. 

Now we prove the assertion we mentioned at the beginning of this 
 subsection.
\begin{lem}
Let $\varphi_0\colon \Bbb P^1\to \Bbb P_1^3$
 be a line of class (2).
Then $\varphi_0$ does not have a first order lift. 
\end{lem}  
\proof
First we consider a line of class (2)-II.
For generic $f$ (the defining polynomial of the quintic hypersurface), 
 as in the calculation in the previous section, 
 the line is given by the equations of the form
\[
x-aw +O_2= 0, \;\; y-bw +O'_2= 0,\;\; z = t= 0,
\]
 here $O_2$ and $O_2'$ are the sum of the terms which are at least 
 quadratic with respect to $x, y, w$,
 and $a, b$ are general nonzero complex constants. 

Take the parameter on $\Bbb P^1$ 
 defined by the pull-back of the function
 $w$.
We write it by $S = \varphi_0^*(w)$.
Then around the point $p = \varphi_0(\Bbb P^1)\cap \ell_{12}$, 
 the pull-back of $x$ and $y$ can be written in the form
\[
\varphi_0^*x = aS + O(S^2),\;\; \varphi_0^*y = bS+O(S^2).
\]

Now assume $\varphi_0$ has a first order lift $\varphi_1$.
Then the pull-back of the functions $w, x, y$ should have 
 the following form:
\[
\varphi_1^*w = S + th_1(S),\;\; \varphi_1^*x = aS + O(S^2)+
 th_2(S),\;\;
 \varphi_1^*y = bS + O(S^2) + th_3(S),
\]
 here $h_i$ are convergent series in the variable $S$.
Similarly, the pull-back of the function $z$ should have the form
\[
\varphi_1^*z = th_4(S),
\]
 here again $h_4$ is a series in $S$.
 
Pulling back the relation $xyz + tw = 0$ by $\varphi_0$,
 the above functions 
 in $S$ must satisfy the relation of the form 
\[
th_4(S)(aS + O(S^2)+ th_2(S))(bS + O(S^2) + th_3(S))
 + t(S + th_1(S)) = 0.
\] 
However, the coefficient of $t$ has the form
\[
S + O(S^2),
\]
 which cannot vanish.
Thus, the map $\varphi_0$ does not have a lift even locally around $p$.

Since the lines of class (2)-I have the same local structure
 at the intersection with $\ell_{j_1j_2}$, this argument shows that
 these curves do not lift, either.\qed\\

\subsubsection{Liftability of lines of class {\rm (1)}}
Here, we prove the following.
It completes the proof of Theorem \ref{thm:CY3}.
\begin{lem}\label{lem:finite}
For generic $f$, the lines of class (1) lift uniquely.
\end{lem}
\proof
Let $\varphi_0\colon \Bbb P^1\to \Bbb P_1^3$ be a line of class (1).
It intersects each of the toric divisors of $\Bbb P_1^3$ at 1 point, 
 and these points must be contained in the singular set $\mathcal S$
 of the total space $\mathfrak X$.

Again by mimicking the calculation in 
 the proof of Lemma \ref{lem:coord}, 
 one sees the following.
Let $p$ be an intersection between $\varphi_0(\Bbb P^1)$ and 
 a toric divisor of $\Bbb P_1^3$. 
Then in a 
 neighborhood $U_p$ of $p$, the total space $\mathfrak X$ is analytically
 isomorphic to a product of  
 a neighborhood of the origin of the set
\[
\{(x, y, z, t)\in\Bbb C^4\;|\; xy+tz = 0\}\subset \Bbb C^3\times\Bbb C
\] 
 and an open disk in $\Bbb C$.
This is the product of the singularity we discussed in the 
 previous section and a disk, 
 and has a natural log structure (as in the previous section, 
 we do not need to define a log structure on the whole $\mathfrak X$, 
 but only in a neighborhood of $\varphi_0(\Bbb P^1)$.
Then the local log structures defined by the local isomorphism to 
 a toric variety as above at each intersection between 
 $\varphi_0(\Bbb P^1)$ and the toric divisors of $\Bbb P_1^3$
 can be extended).
 
So we can apply the calculation 
 in the previous section to study the log normal sheaf of $\varphi_0$.
Namely, let $\mathcal Q$ be the usual
 normal sheaf of the map $\varphi_0$ as a map to 
 $\Bbb P^3\cong\Bbb P_1^3$.
Thus, we have an isomorphism $\mathcal Q\cong \mathcal O(1)\oplus
 \mathcal O(1)$.
Then by the calculation in the previous section (see Lemma \ref{lem:N}),
 the set of local lifts of $\varphi_0$ on the
 neighborhood $U_p$ of $p$ is a torsor over the set of the sections of  
 a sheaf $\mathcal N_{U_p}$ which satisfies the exact sequence
\[
0\to \mathcal N_{U_p}\to \mathcal Q_{U_p}\to \Bbb C_{\varphi^{-1}_0(p)}
 \to 0, 
\] 
 where $\Bbb C_{\varphi^{-1}_0(p)}$ is 
 the skyscraper sheaf on $\Bbb P^1$.

Since we have 4 such intersection points, we see that
 the set of first order lifts of $\varphi_0$ is 
 a torsor over the space of sections of a sheaf $\mathcal N$
 which satisfies the exact sequence
\[
0\to\mathcal N\to \mathcal Q\to \oplus_{i=2}^5\Bbb C_{\varphi^{-1}_0(p_i)}
 \to 0,
\]
 here $p_i$, $i = 2, 3, 4, 5$ are the intersection of $\varphi_0(\Bbb P^1)$
 and the toric divisors of $\Bbb P_1^3$.

\begin{lem}
The sheaf $\mathcal N$ is isomorphic to $\mathcal O(-1)\oplus
 \mathcal O(-1)$.
\end{lem}
\proof
By the exact sequence above, 
 the sheaf $\mathcal N$ is torsion free and has degree $-2$.
Thus, it is isomorphic to either of
\[
\mathcal O(-3)\oplus \mathcal O(1),\;\;
 \mathcal O(-2)\oplus \mathcal O,\;\; \mathcal O(-1)\oplus\mathcal O(-1).
\]
Since the sheaves 
 $\mathcal O(-3)\oplus \mathcal O(1)$ and 
 $\mathcal O(-2)\oplus \mathcal O$ have nontrivial global sections, 
 it suffices to prove that the sheaf $\mathcal N$ does not have
 a nontrivial global section.

Now we can take a hyperplane $H$ of $\Bbb P^3_1$ so that it does not
 intersect the set $\{p_1, p_2, p_3, p_4\}$.
Then take the complement $\Bbb P_1^3\setminus H\cong\Bbb C^3$. 
We can trivialize the usual normal sheaf $\mathcal Q$ of $\varphi_0$
 on $\varphi_0(\Bbb P^1)\cap \Bbb C^3$:
\[
\mathcal Q|_{\Bbb C^3}\cong \mathcal O_{\Bbb P^1\setminus 
 \varphi_0^{-1}(p_0)}\mathfrak v_1
 \oplus \mathcal O_{\Bbb P^1\setminus \varphi_0^{-1}(p_0)}\mathfrak v_2, 
\]
 here $p_0$ is the intersection $\varphi_0(\Bbb P^1)\cap H$
 and $\mathfrak v_i$ are generators of the quotient space
 $\Bbb C^3/\overline{\varphi_0(\Bbb P^1)}$, 
 where $\overline{\varphi_0(\Bbb P^1)}$ is the subspace of 
 $\Bbb C^3$ parallel to the line $\varphi_0(\Bbb P^1)$.

At the points $p_2, p_3, p_4, p_5$, the line 
 $\varphi_0(\Bbb P^1)$ intersects quintic curves 
 $C_2, C_3, C_4, C_5$.
Let $T_i\subset \Bbb C^3/\overline{\varphi_0(\Bbb P^1)}$
 be the subspace spanned by the tangent vector of 
 $C_i$ at $p_i$.
As in the proof of Proposition \ref{prop:transv}, for general $f$, the 
 spaces $T_i$ are all different in 
 $\Bbb C^3/\overline{\varphi_0(\Bbb P^1)}$.
 
With this preparation, now we study the space of sections of $\mathcal N$.
Namely, the restriction of a global section of $\mathcal Q$ 
 to $\Bbb P^1\setminus \varphi_0^{-1}(p_0)$ is described as
\[
f_1(S)\mathfrak v_1+ f_2(S)\mathfrak v_2,
\]
 here $S$ is a parameter on 
 $\Bbb P^1\setminus \varphi_0^{-1}(p_0)$
 and $f_1, f_2$ are affine linear functions.
Then $f_1(S)\mathfrak v_1+ f_2(S)\mathfrak v_2$
 belongs to $\mathcal N$ if and only if
 the condition
\begin{itemize}
\item At each $p_i$, the vector
 $f_1(S)\mathfrak v_1+ f_2(S)\mathfrak v_2$
 is contained in the subspace
 $T_i$
\end{itemize}
 is satisfied.
However, if $f$ is general so that the subspaces
 $T_i$ are also general, by an elementary calculation one sees that 
 such affine functions $f_1$ and $f_2$ do not exist.
Thus, the sheaf $\mathcal N$ does not have a nontrivial global section, 
 and consequently it is isomorphic to 
 $\mathcal O(-1)\oplus \mathcal O(-1)$.\qed\\

Now we finish the proof of Lemma \ref{lem:finite}.
As we mentioned above, the set of first order local lifts of $\varphi_0$ is 
 a torsor over the set of local sections of 
 $\mathcal N = \mathcal O(-1)\oplus \mathcal O(-1)$.
Then the obstruction
\[
H^1(\Bbb P^1, \mathcal O(-1)\oplus \mathcal O(-1))
 \cong H^0(\Bbb P^1, \mathcal O(-1)\oplus \mathcal O(-1))^{\vee}
\]
 vanishes.
Since local lifts of $\varphi_0$ clearly exist, it follows that
 there is a global first order lift of $\varphi_0$ and it is unique.
Similarly, higher order lifts uniquely exists.
This finishes the proof of Lemma \ref{lem:finite}, 
 hence of Theorem \ref{thm:CY3}.\qed

\begin{rem}
The degeneration need not be the ones whose central fiber is a union of 
 projective spaces, but the components of the central fiber 
 can be more general toric variety (in fact, the components of the central fiber 
 need not even
 be toric varieties, though explicit calculations will become difficult in general).
We give examples in the cases of a cubic surface and a 
 quartic surface.\\
 
\noindent
1. Cubic surface.\\

Consider a degeneration of a cubic surface in $\Bbb P^3$
 to the union of a projective plane 
 $P$ and a quadratic surface $Q$.
The intersection $P\cap Q$ is a 
 quadratic curve $C$, and the singular locus $\mathcal S$ of the 
 total space of degeneration is a set of six points on $C$.
There are ${}_6C_2 = 15$ lines on $P$ which intersects 
 $C$ at $\mathcal S$.
On the other hand, recalling $Q\cong \Bbb P^1\times \Bbb P^1$, 
 generically
 there are 12 lines on $Q$ which intersects $C$ at $\mathcal S$
 (there are 2 lines, horizontal and vertical, through each of the 
 six points from $\mathcal S$).
By the same calculation
 as in Section \ref{sec:K3}, one sees that
 the normal sheaf is isomorphic to $\mathcal O(-1)$ in each case, 
 and so all these lines lift, giving 27 lines, as expected.\\
 
\noindent
2. Quartic surface.\\

Consider a degeneration 
 of a quartic surface in $\Bbb P^3$ 
 into the union of 2 quadratic surfaces $Q_1, Q_2$.
The intersection $Q_1\cap Q_2$ is a curve 
 $C$ of degree $(2, 2)$ on each of them.
The singular locus 
 $\mathcal S$ of the total space of degeneration is a set of 
 16 points on $C$.
A line in $Q_i$ intersects $C$ at 2 points (counted with multiplicity), 
 and the necessary condition for the liftability is that both of these intersections
 are contained in $\mathcal S$.
The sufficient condition is given by the vanishing of the Kuranishi map
 which can be calculated as in Subsection \ref{subsec:b_1}.\\

\end{rem}

\section{Disks of Maslov index zero}\label{sec:disks}
Using the technique developed in \cite{N3}, 
 we can apply the calculation in this paper 
 to construct families of holomorphic disks of Maslov index zero in many 
 symplectic manifolds.
Such disks play important roles.
For example: 
\begin{enumerate}
\item When we calculate the \emph{potential functions} (see \cite{FOOO}, Section 4),
 it is important to count disks of Maslov index 2.
These disks undergo wall crossing phenomena, and disks of Maslov index zero
 compose these walls (\cite{A}).
\item Same kind of walls appear in the general construction of Calabi-Yau manifolds
 (\cite{GS, KS}) from affine manifolds with singularities. 
\end{enumerate}
Let $n$ be an integer with $n\geq 3$.
Let $f = f(x_0, \dots, x_n)$ be a general
 homogeneous polynomial of degree larger than 1.
Let $V = \{f = 0\}\subset \Bbb P^n$ be a reduced projective hypersurface.
Let $d$ be the degree of $f$.
Then consider a degeneration
\[
a_0a_1\cdots a_{d-1} + tf  =0
\]
 of $V$.
Here  
 $t$ is 
 the parameter of the degeneration.
We take $a_i$ so that
\[
a_i = \begin{cases}
x_{i},\;\; i = 0, \dots, n.\\
\text{A general linear homogeneous polynomial},\;\; i>n.
\end{cases}
\]
Let $V_t$ be the hypersurface defined by the above equation.
\begin{thm}\label{thm:disk}
Fix a complex number $\varepsilon$
 with $|\varepsilon|$ sufficiently small.
Let 
 $\Omega\in H^2(V_{\varepsilon}, \Bbb Z)$
 be the class determined by the polarization.
Then there is a symplectic form $\omega$
 on $V_{\varepsilon}$ of class $\Omega$
 such that there is a Lagrangian torus $L$ with respect to $\omega$
 which bounds a holomorphic disk of Maslov index zero.
\end{thm}
\proof
Let $\mathfrak X\subset \Bbb P^n\times\Bbb C$
 be the total space of the degeneration
 and
 let
\[
\pi\colon \mathfrak X\to\Bbb C
\]
 be the projection.
We put a symplectic structure on each fiber of $\pi$
 by restriction of the Fubini-Study symplectic form on $\Bbb P^n$
 (on the central fiber $X_0$,
 it suffices to put a symplectic structure on its smooth part).

Let $\mathcal S$ be the singular locus of $\mathfrak X$.
As before, we can assume it is contained in the central fiber $X_0$
 for generic $f$.
Its intersection with a plane 
\[
\Bbb P_{ij}^{n-2} = \Bbb P_i^{n-1}\cap \Bbb P_j^{n-1},\;\;
  i\neq j, \;\; 0\leq i, j\leq n+1
\]
 here $\Bbb P_i^{n-1} = \{a_i = 0, t = 0 \}$,  is a 
 hypersurface of degree $d$ in $\Bbb P_{ij}^{n-2}$.
Let 
\[
\Bbb P_{klm}^{n-3} = \Bbb P_k^{n-1}\cap \Bbb P_l^{n-1}\cap\Bbb P_m^{n-1},\;\; 
 k\neq l\neq m\neq k,\;\; 0\leq k, l, m\leq d
\]
 be the triple intersection.
Let 
\[
int\Bbb P_{ij}^{n-2} = \Bbb P_{ij}^{n-2}\setminus\bigcup_{k\neq i, j}\Bbb P_{ijk}^{n-3}
\]
 be the 'interior' of $\Bbb P_{ij}^{n-2}$.

Let us fix a point $p\in \mathcal S\cap int\Bbb P_{ij}^{n-2}$ generally.
There is a natural structure of a toric variety 
 on $\Bbb P_i^{n-1}$, $0\leq i\leq n$,
 so that $\Bbb P_{ij}^{n-2}$, $0\leq i, j\leq n$, is one of the toric divisors. 
Then by the description of the disks in toric varieties (\cite{CO}), 
 there is a real 1 dimensional family of
 Lagrangian tori in $\Bbb P_i^{n-1}$
 each member of which bounds a unique (up to isomorphism)
 disk of Maslov index 2 
 intersecting the point $p$.

Let $L$ be one of these Lagrangian submanifolds.
It is an orbit of the torus action of $T^{n-1}(\subset (\Bbb C^*)^{n-1})$
 on $\Bbb P_i^{n-1}$. 
The restriction of the 
 Fubini-Study form $\tilde{\omega}_{FS}$ of $\Bbb P^n$ to 
 $\Bbb P_i^{n-1}$ defines a $T^{n-1}$-invariant symplectic form $\omega$.
We want to extend this Lagrangian boundary condition to 
 a family of boundary conditions over the base space, 
 but due to the real nature of the Lagrangian submanifolds, we
 restrict the base space from the complex plane to the real half line
 of appropriate direction so that it contains the point
 $\varepsilon\in\Bbb C$.
Thus, let
\[
\pi_{\Bbb R}\colon \mathfrak X|_{\Bbb R_{\geq 0}}\to \Bbb R_{\geq 0}
\]
 be the restriction of $\pi\colon\mathfrak X\to \Bbb C$.
 
\begin{lem}\label{lem:boundary}
We can take a neighborhood $U$ of $L$ 
 in the total space $\mathfrak X|_{\Bbb R_{\geq 0}}$
 so that the following conditions hold.
\begin{itemize}
\item The diagram
\[\xymatrix{
U \ar[r]\ar[d] & \mathfrak X|_{\Bbb R_{\geq 0}}\ar[d]\\
\Bbb R_{\geq 0} \ar@{=}[r] & \Bbb R_{\geq 0}
}
\]
 commutes, where the upper arrow is the inclusion.
\item The action of $T^{n-1}$ extends to $U$ over $\Bbb R_{\geq 0}$,
 and its orbits are Lagrangian tori with respect to the restriction
 of the Fubini-Study symplectic form on $\Bbb P^n$.
\item On $U$, 
 the action and the form are (real) 
 analytic with respect to the analytic structure of 
 $\mathfrak X$.
\end{itemize}
\end{lem}
\proof
We can take such $U$ as follows.
Let $U'$ be a neighborhood of $L$ in $\Bbb P_i^{n-1}$
 which is invariant under the action of $T^{n-1}$.
Let $[0, \delta)$ be an interval on the base space $\Bbb R_{\geq 0}$
 of the degeneration.
We have a gradient-Hamiltonian flow (\cite{R}, see also \cite{NNU})
\[
\phi_{s}\colon X_{s}\to X_0,\;\; s\in [0, \delta)
\]
 which is defined on the smooth locus of $\mathfrak X$.
The map $\phi_s$ is a diffeomorphism on a dense open subset
 $\phi_s^{-1}(int X_0)$, where $int X_0$ is the complement in $X_0$
 of the union of all toric divisors of the components of $X_0$.
In particular, it is a diffeomorphism in a neighborhood of $L$
 in $\mathfrak X$. 
Moreover, it
 preserves the symplectic structures on the fibers.
 
Then define $U$ by
\[
U = \cup_{0\leq s <\varepsilon}\phi_s^{-1}(U').
\]
Put the action of $T^{n-1}$ on $U$ by pulling back the action of it 
 on $U'$.
Since $\phi_s$ preserves the fiberwise symplectic structure, 
 the orbits of $T^{n-1}$ are Lagrangian submanifolds
 of $X_s$.
The analyticity of the action and the symplectic forms are clear since 
 the gradient-Hamiltonian flow is real analytic.\qed\\

Thus, we have a family $U$ of Lagrangian tori containing $L$
 each of which can be 
 a boundary condition for holomorphic disks.
In \cite{N3}, we developed a deformation theory of disks with 
 such Lagrangian boundary conditions.
Let 
\[
\varphi_0\colon D\to X_0
\]
 be the disk of Maslov index 2 with boundary on $L$ intersecting
 the point $p\in\mathcal S\cap \Bbb P_{ij}^{n-1}$ defined above.
Then by \cite{N3}, the lifts of $\varphi_0$ to the fibers over $t\neq 0$ are
 controlled by the log normal sheaf
 as in the case of rational curves,
 except that we should use 
 sheaves associated to \emph{Riemann-Hilbert bundles} (\cite{KL})
 to take care of the boundary condition.

In our case, the log normal bundle of the disk $\varphi_0$, 
 as a map to $\Bbb P_i^{n-1}$, is trivial as a Rimenn-Hilbert bundle.
However, as in the case of rational curves, we have to take account of the 
 singularity of $\mathfrak X$, which contains the point $p$.

For simplicity, let us take $i = 0$, $j = 1$, and use the affine coordinates
\[
Y_i = x_i/x_n,\;\; i = 0, \dots, n-1
\]
 on $\Bbb P^n$.
We also set
\[
Y_j = a_j/x_n
\]
 for $j\geq n+1$.
Using the coordinates $\{Y_i\}_{i=0}^{n-1}$, the point $p$ is represented as
\[
p = (0, 0, b_2, \dots, b_{n-1}),
\]
 where $b_i\neq 0$ since $p$ lies in $int\Bbb P_{01}^{n-2}$.
Also, the linear polynomials $Y_i$, $n+1\leq i\leq d$ are not zero at $p$.
Note that the disk $\varphi_0(D)$ 
 mentioned above is parametrized
 by
\[
(s, 0, b_2, \dots, b_{n-1}),
\]
 where $s$ is an affine coordinate on the disk.
Redefining $f$ by 
\[
\frac{f}{Y_2Y_3\cdots Y_d},
\]
 the defining equation of the degeneration
 near $p$ is
\[
Y_0Y_1 + tf = 0.
\]
Here the function $f$ is also written in the affine coordinates 
 $Y_i$, $i = 0, \dots, n-1$.

Since $p$ belongs to the singular locus, 
\[
f(p) = 0
\]
holds.
By genericity of $f$ and $p$, we can assume the following:
\begin{itemize}
\item 
$
\frac{\partial f}{\partial Y_0}(p)\neq 0.
$
\item For at least one of $i\in \{2, 3, \dots, n-1\}$, 
$
\frac{\partial f}{\partial Y_i}(p)\neq 0.
$
\end{itemize}
Under this assumption, the function $f$ is pulled back 
 by $\varphi_0$ to a holomorphic function on the disk which
 has a simple zero at $\varphi^{-1}_0(p)$.
Also, we can use the functions
\[
Y_0, Y_1, \dots, Y_{i-1}, f, Y_{i+1}, \dots, Y_{n-1}, t
\]
 as a coordinate system of $\mathfrak X$ in a neighborhood of $p$.
In particular, in a neighborhood of $p$, $\mathfrak X$ is 
 analytically isomorphic to the product of 
 a neighborhood of the origin of 
\[
\{(X, Y, Z, t)\in\Bbb C^4\;|\; XY+tZ = 0\}\times O
 \subset (\Bbb C^3\times\Bbb C)\times O,
\]
 where $O$ is an open neighborhood of the origin in the affine space
 with a coordinate system 
 $\{Y_2, Y_3, \dots, Y_{i-1}, Y_{i+1}, \dots, Y_{n-1}\}$.

Differentiating the equation $Y_0Y_1+tf = 0$, we obtain the relation
\[
\frac{dY_0}{Y_0}+\frac{dY_1}{Y_1}-\frac{df}{f}-\frac{dt}{t} = 0.
\]
Restricting to $X_0$, we have
\[
\frac{dY_0}{Y_0}+\frac{dY_1}{Y_1}-\frac{df}{f} = 0.
\]
We separate the problem into 2 
 cases according to $n = 3$ or $n\geq 4$.
When $n$ equals to 2, the set $\mathcal S$ is empty and there is no
 disk of Maslov index 0.\\

\noindent
{\bf Case 1: $n = 3$.}
In this case, the singular locus of $\mathfrak X$ is zero dimensional, and 
 the argument in Section \ref{sec:K3} extends straightforwardly.
Let us choose an analytic family of Lagrangian tori
 $\mathfrak L\subset \mathfrak X$ over $\Bbb R_{\geq 0}$
 such that
\begin{itemize}
\item on $X_0$, it restricts to $L$ fixed above, 
 and 
\item each fiber over $t\in\Bbb R_{\geq 0}$
 is an orbit of the $T^{n-1}$-action given in Lemma \ref{lem:boundary}.
\end{itemize} 
We use $\mathfrak L$ as the boundary condition for the deformation
 of $\varphi_0$.
Then we have the following.
\begin{lem}\label{lem:discnormal}
The set of first order local lifts of the map $\varphi_0$ is a torsor over 
 the set of local sections of  
 the sheaf 
 associated to a Riemann-Hilbert bundle whose doubling is
 isomorphic to $\mathcal O(-2)$.
\end{lem} 
\proof
The usual normal sheaf as a map to 
 $\Bbb P_0^{n-1}$ is the sheaf 
 $(\varphi_0^*\Theta_{X_0}, \varphi_0^*\Theta_L)/
    (\Theta_D, \Theta_{\partial D})$ associated to 
 the quotient Riemann-Hilbert bundle,
 where $\Theta_L$ is the sheaf of real analytic sections of 
 the tangent bundle of $L$,
 $\Theta_D$ is the restriction to $D$ of the sheaf of sections
 of the tangent bundle of the complex plane $\Bbb C$, 
 and $\Theta_{\partial D}$ is the sheaf of real analytic sections of 
 the tangent bundle of the boundary $\partial D$ of $D$, 
 see \cite{N3}.
This is a sheaf associated to a trivial Riemann-Hilbert bundle of 
 rank $1$.

As we calculated in Lemma \ref{lem:N}, 
 the set of local lifts of $\varphi_0$ on a neighborhood of the origin of $D$,
 which is mapped to the singular locus of $\mathfrak X$, 
 is a torsor over the set of local sections of the usual normal sheaf which
 vanish at the origin.
This contributes $-1$ to the Riemann-Hilbert bundle,
 and then it is doubled when we take the double.\qed\\
 
\begin{defn}
We write by $\mathcal N$ the sheaf mentioned in
  Lemma \ref{lem:discnormal}.
Namely, the set of sections of $\mathcal N$ is the same as the set
 of sections of  
 the sheaf
 $(\varphi_0^*\Theta_{X_0}, \varphi_0^*\Theta_L)/
    (\Theta_D, \Theta_{\partial D})$,
 with the extra condition that it vanishes at the origin of $D$.
\end{defn}

Let $\omega$ be the sheaf 
 associated to the canonical Riemann-Hilbert bundle of the disk, 
 whose doubling is also $\mathcal O(-2)$.
By Serre duality, we have
\[
H^1(C_0, \mathcal N)
 \cong H^0(C_0, \mathcal N^{\vee}\otimes \omega)^{\vee}
 \cong H^0(C_0, (\mathcal O_D, \mathcal O_{\partial D}))^{\vee}
 \cong \Bbb R,
\]
 here $\mathcal O_D$ is 
 the restriction of $\mathcal O_{\Bbb P^1}$ to $D$ and
 $\mathcal O_{\partial D}$ is the real analytic subsheaf of 
 $\mathcal O_D|_{\partial D}$ composed of real valued functions
 (see \cite{N3}, Section 7).

By the same calculation as in Section \ref{sec:K3}, we can 
 describe 
 the generator of $H^0(C_0, \mathcal N^{\vee}\otimes \omega)$
 explicitly.
Namely, it is given by a logarithmic form
\[
\eta\otimes\frac{ds}{s},
\]
 where $s$ is a coordinate on $D$ such that the point $s = 0$
 is mapped to the singular locus of $\mathfrak X$,
 and $\eta$ is a constant section of the usual conormal bundle of
 $\varphi_0$, which is trivial.

The obstruction to lift $\varphi_0$ is calculated as in
 Theorem \ref{thm:Kuranishi}.
Namely, take an open covering $\{O_i\}$
 of the domain $D$ of $\varphi_0$
 so that on each open subset a local lift of $\varphi_0$ 
 satisfying the boundary condition $\mathfrak L$ exists. 
Such a local lift of $\varphi_0$ 
 corresponds to a local section of the normal sheaf of $\varphi_0$
 as a map to $\mathfrak X$.
Explicitly, those local sections can be written in the form
 analogous to
 $(\ast)$ in the proof of Theorem \ref{thm:Kuranishi}.
 
Then take the pairing 
 between that section and the $\mathcal N^{\vee}$-part of
 the form $\eta\otimes\frac{ds}{s}$ so that we have a locally defined
 1-form.
In the present case, these 1-forms have nontrivial residue
 only at the point $s = 0$, 
 and the value of the residue does not depend on the lift.

The difference from the argument in Section \ref{sec:K3}
 is, in addition to such a residue, 
 we have to take care of the boundary.
However, the contribution from the boundary is calculated in the
 same manner.
Namely, take the pairing of the local lift with the form 
 $\eta\otimes\frac{ds}{s}$ as above,
 then integrate the resulting 1-form along the boundary
 (multiplied by $\frac{1}{2\pi i}$).

To be more precise, 
 let us triangulate the disk $D$ so that
\begin{itemize}
\item the origin of $D$ is contained in the interior of some triangle, and
\item each triangle is contained in an open subset $O_i$
 of $D$ introduced above.
There may be several $O_i$ containing a given triangle, then 
 choose one of them.
We call this $O_i$ the open subset associated to the given triangle. 
\end{itemize}
Also, we fix a direct sum decomposition of 
 $\varphi_0^*\Theta_{\mathfrak X}$ into base and fiber direction:
\[
\varphi_0^*\Theta_{\mathfrak X} = \mathcal O_{C_0}\cdot t\partial_t
 \oplus \Theta_{\mathfrak X}^f.
\]
Correspondingly, the normal sheaf 
 $\varphi_0^*\Theta_{\mathfrak X}/\Theta_{C_0}$
 of $\varphi_0$ as a map to 
 $\mathfrak X$ decomposes as
\[
\varphi_0^*\Theta_{\mathfrak X}/\Theta_{C_0}
 = \mathcal O_{C_0}\cdot t\partial_t\oplus \mathcal N.
\]
Then we have the following.
\begin{lem}\label{lem:pairing}
The pairing between the obstruction to lift $\varphi_0$ 
 $($which lies in $H^1(C_0, \mathcal N))$ and the class of 
 $H^0(C_0, \mathcal N^{\vee}\otimes \omega) (\cong
 H^1(C_0, \mathcal N)^{\vee})$
 defined by the form $\eta\otimes \frac{ds}{s}$ 
 is calculated by the
 following procedure.
\begin{itemize}
\item Take the pairing between the $\mathcal N$-component of the
 local section of 
 the normal sheaf $\varphi_0^*\Theta_{\mathfrak X}/\Theta_{C_0}$
 on each $O_i$ corresponding to the local lift, 
 and the form $\eta\otimes \frac{ds}{s}$.
Then on each triangle of the above triangulation, 
 we have a 1-form on it by pulling back the 1-form on the associated
 open subset $O_i$.
\item
Integrate the resulting local 1-form along the boundary of the triangle
 which contains the origin of $D$.
Note that the triangle has a natural orientation since it is
 contained in $D\subset \Bbb C$.
\item When a triangle intersects $\partial D$, 
 then integrate the above 1-form along its intersection with $\partial D$
 in the reverse direction of the given orientation.
\item Sum up all the values of the above integrals. 
\end{itemize}
\end{lem}
\proof
It is straightforward to see that the number obtained by the 
 given procedure does not depend on various choices 
 (triangulation, open covering $\{O_i\}$, choice of $O_i$ to 
 each triangle, choices of local lifts, and the decomposition
 $\varphi_0^*\Theta_{\mathfrak X} = \mathcal O_{C_0}\cdot t\partial_t
 \oplus \Theta_{\mathfrak X}^f$
 ).
Thus, this gives a well-defined pairing between the cohomology classes
 $H^1(C_0, \mathcal N))$ and  
 $H^0(C_0, \mathcal N^{\vee}\otimes \omega)$.
Also, it is easy to see that this pairing is not trivial,
 so it is the natural paring between the dual spaces
 possibly up to a constant multiple.
  \qed\\

In our situation, the contribution from the point $s = 0$
 is fixed from the data of the given degeneration $\mathfrak X$, 
 as in the calculation in Section \ref{sec:K3}.
This is given by the terms $b_i$ in the calculation of
 Subsection \ref{subsec:b_1}.
So for the existence of a lift of $\varphi_0$, we have to 
 cancel this by the contribution from the boundary.
This corresponds to moving the boundary Lagrangian torus 
 by the amount determined by $b_i$.

Explicitly, let $X_{0, \alpha}$ be the component of $X_0$
 to which the disk is mapped.
Let $p$ 
 be the point on $X_{0, \alpha}$, which is the intersection
 of the disk with the singular locus of $\mathfrak X$. 
Take a coordinate system $\{x, y\}$ on $X_{0, \alpha}$ around $p$.
We can take $x, y$ so that the following properties are satisfied:
\begin{itemize}
\item The functions $x$ and $y$ are characters of the torus 
 action on $X_{0, \alpha}$.
Note that these are also a part of a coordinate system of
 $\Bbb P^3\supset X_t$.
Thus, we can talk about the value of these functions on the fibers
 $X_t$ of $\mathfrak X$.
\item Let $(x, y) = (0, c)$ be the coordinate of the point $p$.
Then the image of the disk is parametrized as
\[
(s, c),
\]
 where $s$ is the coordinate on the disk given by
\[
s = \varphi_0^*(x),
\]
 see \cite[Proposition 7.3]{CO}.
\end{itemize}
The boundary Lagrangian torus is fixed by
 specifying the norms of the values of $x$ and $y$.
 
In our case, the norm of $y$ of the image of $\varphi_0$
 is fixed to $|c|$.
As in the calculation in Subsection \ref{subsec:b_1}, 
 a local lift of $\varphi_0$ around $s = 0$ gives a term $b_1$, 
 which perturbs the value of the coordinate $y$ from $c$ to $c+tb_1$.
Namely, a local lift has a parametrization of the form
\[
\begin{array}{l}
x = s,\\
y = c+ t(f_1(s)+b_1),
\end{array}
\]
 here $f_1(s)$ is a convergent series without constant term.
Then the Lemma \ref{lem:pairing} shows that
 for the existence of a lift of $\varphi_0$ to a map over $\Bbb C[t]/t^2$, 
 it is enough to take the boundary torus so that 
 the $y$-coordinates
 of points on that torus are $c+tb_1$
 times $e^{i\theta}$, $0\leq \theta<2\pi$.
\begin{rem}
Precisely, the term $b_1$ should be interpreted as 
 the coefficient of the vector
 $b_1\partial_y$, where $\partial_y$ is defined
 using the coordinate system
 $\{x, y\}$.
If the function $y$ is multiplied by a non-zero constant, the constant
 $b_1$ is changed, but the vector $b_1\partial_y$ is not changed.
\end{rem}
\begin{rem}\label{rem:boundaryLag}
When we determine the y-coordinates of points
 on the torus as above,
 the boundary Lagrangian torus is uniquely determined
 when we fix the norm of the $x$-coordinate.
\end{rem}

Generally, if we have a $k$-th order lift $\varphi_k$ of $\varphi_0$, 
 then we can construct a $k+1$-th order lift by taking an appropriate
 boundary condition.
Such a boundary condition can be fixed in the following way.
Namely, take a section $\sigma$ of the projection
 $\pi|_{\Bbb R}\colon \mathfrak X|_{\Bbb R_{\geq 0}}
  \to \Bbb R_{\geq 0}$ restricted to a map over $\Bbb R[t]/t^{k+2}$
 so that 
 the $y$-coordinate of the image of $\sigma$ is given by 
\[
c+tb_1+t^2b_2+\cdots + t^{k+1}b_{k+1},
\]
Here the complex numbers $b_i$ are given by the constant
 terms of a parametrization of a local lift of $\varphi_0$:
\[
\begin{array}{l}
x = s,\\
y = c+t(f_1(s)+b_1)+t^2(f_2(s)+b_2)+\cdots + t^k(f_{k+1}(s)+b_{k+1}).
\end{array}
\]
Here $f_i(s)$ is a convergent series
 without constant terms, and $b_{k+1}$ is determined 
 uniquely when we fix a $k$-th order lift $\varphi_{k}$ of $\varphi_0$.

Then take the family of Lagrangian tori $\mathcal L_{k+1}$
 over $\Bbb R[t]/t^{k+2}$
 given by taking the
 orbit of $\sigma$ by the action of $T^{n-1}$.
We summarize the result so far. 
\begin{prop}\label{prop:n=3}
Take the family of Lagrangian tori $\mathcal L_{k+1}$ as above.
Then there is a map
\[
\varphi_{k+1}\colon D\times Spec\; \Bbb R[t]/t^{k+2}|_{\Bbb R_{\geq 0}}\to 
 \mathfrak X|_{\Bbb R_{\geq 0}}
\]
 over $Spec\; \Bbb R[t]/t^{k+2}$ (restricted to the non-negative part)
 which coincides with $\varphi_k$ when specialized to a map
 over $Spec\;\Bbb R[t]/t^{k+1}$.
Moreover, the map $\varphi_{k+1}$ is unique up to automorphisms. 
\end{prop}
\proof
The existence follows because we constructed $\mathcal L_{k+1}$
 so that the obstruction to lift $\varphi_{k}$ to the $k+1$-th order
 vanishes.
The uniqueness follows from the fact that 
 $H^0(C_0, \mathcal N) = 0$, where
 $\mathcal N$ is the Riemann-Hilbert sheaf introduced above.\qed\\

By a suitable implicit function theorem for analytic category (\cite{Ar}), 
 this proposition shows the existence of a disk of Maslov index zero
 on $V_{\varepsilon}$.
This proves Theorem \ref{thm:disk} for the case $n = 3$.\\

\noindent
{\bf Case 2: $n \geq 4$.}
Note that by a suitable coordinate change as in the previous section, 
 one sees that a neighborhood of a general singular
 point $x\in\mathfrak X$
 is analytically isomorphic to a product of 
 a neighborhood $N$ of the origin of the set
\[
\{(X, Y, Z, t)\in\Bbb C^4\;|\; XY+tZ = 0\}\subset \Bbb C^3\times\Bbb C
\] 
 and an open subset $O$ of $\Bbb C^{n-3}$.
In particular, the pull-back of $\Theta_{\mathfrak X}$ by $\varphi_0$ and 
 the log normal sheaf of $\varphi_0$ defined using the local log structure 
 induced from this toric structure
 are locally free.
Here, we put a log structure (in the analytic category)
 on the disk associated to the point $\varphi_0^{-1}(p)$
 seen as a divisor of degree 1.

Assume that the point $p$ is contained in $\Bbb P_{01}^{n-2}$
 and the map $\varphi_0$ maps the disk $D$ to 
 the component $\Bbb P_0^{n-1}$ of $X_0$.
Then since we take $p$ generic, the homogeneous coordinate
 $x_n$
 is not 0 at $p$.
So we have a standard local coordinate system 
\[
Y_0, Y_1, \dots, Y_{n-1}, t
\]
 of $\Bbb P^n\times \Bbb C$ around $p$, 
 using the notation introduced before.
We write its associated generators of the tangent sheaf by
\[
\partial_{Y_0}, \partial_{Y_1}, \dots, \partial_{Y_{n-1}}, \partial_t.
\]

On the other hand, 
 since $f$ is general, at least one of the derivatives
 $\frac{\partial f}{\partial Y_2}, \dots, \frac{\partial f}{\partial Y_{n-1}}$
 is not 0 around $p$
 (here we abused the notation a little
 and represent by the same letter
 $f$ the function which is the inhomogenized one 
 using the affine coordinates $Y_0, Y_1, \dots, Y_{n-1}$ on $\Bbb P^n$).
Let us assume $\frac{\partial f}{\partial Y_{n-1}}$ is not 0.
Then the functions 
\[
Y_0, Y_1, \dots, Y_{n-2}, f, t
\]
 also give a coordinate system 
 of $\Bbb P^n\times \Bbb C$ around $p$.
This coordinate system gives an explicit isomorphism 
 between a neighborhood of $p$ and a neighborhood of the origin of 
 the model 
 $\{(X, Y, Z, t)\in\Bbb C^4\;|\; XY+tZ = 0\}\times O
 \subset (\Bbb C^3\times\Bbb C)\times O$
 above (the coordinates $Y_2, \dots, Y_{n-2}$
 correspond to the coordinates on 
 the open subset $O$ of $\Bbb C^{n-3}$).

Consider the family of subvarieties 
 $\{S_t\}_{t}$ in $\mathfrak X$ defined by the
 equations
\[
Y_2 = c_2,\;\; Y_3 = c_3, \;\; \cdots, \;\; Y_{n-2} = c_{n-2},
\]
 where $c_i$ are constants.
For each $t$, this gives a subvariety of dimension 2 in $X_t$
 defined in a neighborhood of $p$ in $\mathfrak X$.

Let us write the coordinate of $p$ by
\[
(Y_0, Y_1, Y_2, \cdots, Y_{n-2}, Y_{n-1}, t) = 
 (0, 0, b_2, \cdots, b_{n-2}, b_{n-1}, 0).
\]
Then taking $c_i = b_i$, $i = 2, \dots, n-2$,
 the surface $S_0$ contains the disk $\varphi_0(D)$.
We consider the deformation of $\varphi_0(D)$ in 
 the family $\{S_t\}_t$. 

Since we take $p$ and $f$ generic, the surface $S_0$ is transversal
 to the singular locus $\mathcal S$ at $p$. 
Then, in the family $\{S_t\}_t$,
 the calculation reduces to the case when $n = 3$
 and by Proposition \ref{prop:n=3}, 
 we have a unique lift of $\varphi_0$ when we fix the boundary condition.
This proves the theorem.\qed\\

In fact, in the proof of the theorem 
 we essentially proved more than the 
 statement of the Theorem \ref{thm:disk}. 
Namely, in the case of $n = 3$, by Remark \ref{rem:boundaryLag}, 
 the boundary Lagrangian is uniquely determined when we fix
 the norm of the $Y_0$-coordinate.
This is the same for the family $\{S_t\}$ in the case of $n\geq 4$.
On the other hand, the families of surfaces $\{S_t\}$ are
 parametrized by the constants $c_2, \dots, c_{n-2}$.
On each of these families, which is locally isomorphic to 
 the family $\mathfrak X$ in the case of $n = 3$, 
 the consequence of Proposition \ref{prop:n=3} holds.
Thus, the disks also compose an $n-3$ dimensional family.

Since the disk is of real dimension 2, the total space of these families is
 (real) codimension 2 in $X_t$.
This can be thought of as the geometric realization of the \emph{walls} of 
 \cite{GS}.
We summarize this as follows.
\begin{cor}
The disks constructed in Theorem \ref{thm:disk} compose
 $n-3$ dimensional family.
In particular, the union of these disks composes a submanifold of
 $V_{\varepsilon}$ of real codimension 1.\qed 
\end{cor}


\begin{thebibliography}{99} 
\bibitem{Ar}{\sc Artin, M., }
{\it Algebraization of formal moduli I},
 Global Analysis (Papers in Honor of K. Kodaira), 
 21--71, Univ. Tokyo Press, 1969.
\bibitem{A}{\sc Auroux,D.,} 
{\it Special Lagrangian fibrations, wall-crossing, and mirror symmetry.}
In Surveys in differential geometry. Vol. XIII. Geometry, analysis, and algebraic
geometry: forty years of the Journal of Differential Geometry, volume 13 of Surv. Differ. Geom., pages 1-47. Int. Press, Somerville, MA, 2009.
\bibitem{CO}{\sc Cho,C. and Oh,Y.,}
{\it Floer cohomology and disc instantons of Lagrangian
torus fibers in Fano toric manifolds.}
 Asian J. Math., 10 (2006), no.4, 773-814.
\bibitem{F}{\sc Forster,O.,}
{\it Lectures on Riemann surfaces.}
Graduate Text in Math., Vol. 81. (1999), Springer Verlag.
\bibitem{FOOO}{\sc Fukaya,K., Oh,Y., Ohta,H. and Ono,K.,}
{\it Lagrangian Floer theory on compact toric manifolds I.}
 Duke Math. J. Volume 151, Number 1 (2010), 23-175.
\bibitem{GS}{\sc Gross,M. and Siebert,B.,}
{\it From real affine geometry to complex geometry.}
 Ann. of Math. (2), 174 (2011), no. 3, 1301-1428.
\bibitem{KK}{\sc Kato,K.,}
{\it Logarithmic structures of Fontaine-Illusie.}
Algebraic analysis, geometry, and number theory 
(Baltimore, MD, 1988),  191--224,
 Johns Hopkins Univ. Press, Baltimore, MD, 1989. 
\bibitem{Katz0}{\sc Katz,S.,}
{\it Degenerations of quintic threefolds and their lines.}
Duke Math. J. 50 (1983), no. 4, 1127-1135.
\bibitem{Katz}{\sc Katz,S.,}
{\it Lines on complete intersection threefolds with $K=0$.}
Math. Z. 191 (1986), no. 2, 293-296.
\bibitem{KL}{\sc Katz,S. and Liu,C.,} 
{\it Enumerative geometry of stable maps with
 Lagrangian boundary conditions and multiple covers of the disc.}
 Adv. Theor. Math. Phys.  5  (2001),  no. 1, 1--49. 
\bibitem{KS}{\sc Kontsevich,M. and Soibelman,Y.,}
{\it Affine structures and non-Archimedean analytic spaces.} 
The unity of mathematics (P. Etingof, V. Retakh, I.M. Singer, eds.), 
 321-385, Progr. Math. 244,
Birkh\"auser 2006.
\bibitem{L}{\sc Li,J.,}
{\it Stable morphisms to singular schemes and relative stable morphisms.}
 J. Differential Geom. 57
(2001), 509-578.
\bibitem{Ma}{\sc Manetti,M.,}
{\it Deformation theory via differential graded Lie algebras.}
 arXiv:math/0507284.
\bibitem{N}{\sc Nishinou,T.,}
{\it Correspondence theorems for tropical curves I.}
Preprint.
\bibitem{N2}{\sc Nishinou,T.,}
{\it Toric degenerations, tropical curves and
 Gromov-Witten invariants of Fano manifolds.}
Canad. J. Math. 67 (2015), no. 3, 667-695.
\bibitem{N3}{\sc Nishinou,T.,}
{\it Disk counting on toric varieties via tropical curves.}
Amer. J. Math., 134 (2012), 1423--1472.
\bibitem{N4}{\sc Nishinou,T.,}
{\it Rational curves in Fano hypersurfaces and tropical curves.}
Preprint.
\bibitem{N5}{\sc Nishinou,T.,}
{\it Degeneration and curves on K3 surfaces.}
Preprint.
\bibitem{NNU}{\sc Nishinou,T., Nohara,Y. and Ueda,K.,}
{\it Toric degenerations of Gelfand-Cetlin systems and potential
functions.} Adv. Math. 224 (2010), no. 2, 648--706.
\bibitem{NS}{\sc Nishinou,T. and Siebert,B.,}
{\it Toric degenerations of toric varieties and tropical curves.}
 Duke Math. J. 135 (2006), no. 1, 1--51.
\bibitem{NY}{\sc Nishinou,T. and Yu,T.,}
{\it Realization of tropical curves in abelian surfaces.}
In preparation. 
\bibitem{R}{\sc Ruan, W-D.,} 
 {\it Lagrangian torus fibration of quintic hypersurfaces. {I}.
              {F}ermat quintic case.}
    Winter School on Mirror Symmetry, Vector Bundles and
              Lagrangian Submanifolds (Cambridge, MA, 1999),
    AMS/IP Stud. Adv. Math. vol 23, 297--332.
    Amer. Math. Soc. Providence, RI, 2001.
 
 
 
 
 
 
 
 
 
 
 
 
 
 
 
 
 


\end{thebibliography}
\end{document}